\documentclass[12pt,reqno]{amsart}
\usepackage{amsmath,amsthm,amsfonts,amssymb,times}
\usepackage{verbatim}
\usepackage{url}
\usepackage{latexsym}
\usepackage{mathtools}
\usepackage{graphicx}
\usepackage{url}	

\usepackage[backref = page]{hyperref}
\hypersetup{
    colorlinks=true,
    linkcolor=blue,
    urlcolor =blue, 
    citecolor = magenta}

\renewcommand*{\backrefalt}[4]{%
\ifcase #1 %
{\color{red} No citations.}%
\or
(p.~#2).%
\else
(pp.~#2).%
\fi
}

\usepackage[nameinlink, capitalize, noabbrev]{cleveref}
\crefname{type}{type}{types}

\usepackage[utf8]{inputenc}
\usepackage[T1]{fontenc}
\usepackage{amsfonts,color}
\usepackage[greek,english]{babel}
\setbox0=\hbox{$+$}
\newdimen\plusheight
\plusheight=\ht0
\def\+{\;\lower\plusheight\hbox{$+$}\;}

\setbox0=\hbox{$-$}
\newdimen\minusheight
\minusheight=\ht0
\def\-{\;\lower\minusheight\hbox{$-$}\;}

\setbox0=\hbox{$\cdots$}
\newdimen\cdotsheight
\cdotsheight=\plusheight
\def\cds{\lower\cdotsheight\hbox{$\cdots$}}

\makeatletter
\def\leqalignno#1{\displ@y \tabskip\z@ plus\@ne fil
  \halign to\displaywidth{\hfil$\@lign\displaystyle{##}$\tabskip\z@skip
    &$\@lign\displaystyle{{}##}$\hfil\tabskip\z@ plus\@ne fil
    &\kern-\displaywidth\rlap{$\@lign\hbox{\rm##}$}\tabskip\displaywidth\crcr
    #1\crcr}}
\makeatother

\renewcommand{\Im}{\operatorname{Im}}
\newcommand{\df}{\dfrac}

\newcommand\pieter[1]{{\textcolor{black}{#1}}}

\renewcommand{\Re}{\text{Re}}
\renewcommand{\Im}{\text{Im}}
\renewcommand{\(}{\left\(}
\renewcommand{\)}{\right\)}
\renewcommand{\[}{\left\[}
\renewcommand{\]}{\right\]}
\numberwithin{equation}{section}
\theoremstyle{plain}
\newtheorem{theorem}{Theorem}[section]

\newtheorem{corollary}[theorem]{Corollary}

\newtheorem{entry}[theorem]{Entry}
\newtheorem{conjecture}[theorem]{Conjecture}
\newtheorem{remark}[theorem]{Remark}
\newtheorem{definition}[theorem]{Definition}
\newtheorem{example}[theorem]{Example}
\newtheorem{claim}[theorem]{Claim}

\makeatletter
\newcommand{\leqnomode}{\tagsleft@true\let\veqno\@@leqno}
\newcommand{\reqnomode}{\tagsleft@false\let\veqno\@@eqno}
\makeatother
\renewcommand{\pmod}[1]{\,(\textup{mod}\,#1)}
\allowdisplaybreaks

\begin{document}
\title{Sums of two squares and the tau-function: Ramanujan's trail}
\author{Bruce C.~Berndt and Pieter Moree}
\address{Department of Mathematics, University of Illinois at Urbana-Champaign, Urbana, Illinois 61802, U.S.A.}
\email{berndt@illinois.edu}
\address{Max Planck Institute for Mathematics, Vivatsgasse 7, 53111, Bonn, Germany}
\email{moree@mpim-bonn.mpg.de}

\subjclass{01A99, 11N37, 11F33}
\keywords{Ramanujan, sums of two squares, modular forms, divisibility of Fourier coefficients, }
\date{}

\begin{abstract}
\noindent
Ramanujan, in his famous first letter to 
Hardy, claimed a
very precise estimate for the number of integers that can be written as
a sum of two squares. Far less well-known is that
he also made further claims of a similar nature for the non-divisibility
of the Ramanujan tau-function
for
certain primes. In this survey, we provide more historical details and also discuss related later developments. These show that, as
so often, Ramanujan was an explorer in a fascinating
wilderness, leaving behind him a beckoning trail.
\end{abstract}

\maketitle

\vspace{-0.5cm}
\section{Ramanujan and sums of two squares}

\label{Ramasquare}
\noindent G.H.~Hardy called his collaboration with Ramanujan ``the one romantic incident
of my life'' \cite[p.~2]{Hardy}.  The foundation for this collaboration
was laid in 1913 when Hardy received  a letter from Ramanujan,\footnote{
The complete correspondence between Ramanujan and Hardy can be found in the book by the first author and Robert Rankin
\cite[Chapter~2]{Letters}. For biographic information about Ramanujan, Kanigel's book \cite{Kanigel} is our favorite.} in which the latter stated  many
results that he claimed to have proved.  Of some of them, Hardy remarked that they
could only have been written down by a mathematician of the highest class,
and moreover that ``They must be true because, if they were not true, no one
would have had the imagination to invent them'' \cite[p.~9]{Hardy}. The claim in Ramanujan's first letter that is
our entry point into Ramanujan's trail reads:

\begin{claim}  \label{ramaclaim}
The number of numbers between $A$ and $x$ which are either squares or
sums of two squares is
$$K\int_{A}^x \frac{dt}{\sqrt{\log t}}+\theta(x),$$
where $K=0.764\dotsc$ and $\theta(x)$ is very small compared with
the previous integral.  $K$ and $\theta(x)$ have been exactly found, though
complicated \dots
\end{claim}
Answering an inquiry of Hardy,
Ramanujan wrote in his second
letter \cite[p.~56]{Letters}: ``The order of $\theta(x)$ which you asked in your letter is
$\sqrt{x/\log x}$.'' 
In his book \emph{Ramanujan} \cite[p.~61]{Hardy}, Hardy informs us that Ramanujan also gave the exact value of $K$, namely,
\begin{equation}\label{constantK}
K=\frac{1}{\sqrt2}\prod_{p\equiv 3\pmod{4}}\Bigl(1-\frac{1}{p^2}\Bigr)^{-1/2}.
\end{equation}
(Here and in the remainder of the paper the letter $p$ is exclusively used to indicate a prime number.)
Using Euler's evaluation of $\zeta(2)$ and its
product representation (now called \emph{Euler product}) for it, i.e.,
$$\zeta(2)=\sum_{n\ge 1}\frac{1}{n^2}=\prod_p\frac{1}{1-p^{-2}}=\frac{\pi^2}{6},$$
we see that $K$ can alternatively be written as
$$
    K=\frac{\pi}{4}\prod_{p\equiv 1\pmod{4}}\Bigl(1-\frac{1}{p^2}\Bigr)^{1/2}.
$$
Ramanujan recorded the value of $K$ twice in his notebooks \cite[Vol.~2, pp.~307, 363]{nb}.  In the first edition of \cite{nb}, the entry on page 307 is very difficult to read.  Moreover, in addition to its lack of clarity in the first edition, the value $K$ ($C$ in his notation on page 307) is incomplete.
 But in the second edition, the entry is much clearer; in the pagination of the second edition, the entry appears on page 350. In both editions, the entry on page 363 is clear.

Both Hardy and G.N.~Watson expressed their
wonder as
to how Ramanujan, at such a young age and isolated from the centers of mathematics, came  to
discover such a
result, which seems to require knowledge of complex analysis for its resolution,
which he very likely did not have at that time.
In particular, Hardy asserted \cite[p.~xxv]{cp}:
\begin{quote}
The dominant term, viz.~$KB(\log\,B)^{-1/2}$, in Ramanujan's notation, was first obtained by 
E.~Landau in 1908.  Ramanujan had none of Landau's weapons at his command; \dots. It is sufficiently marvellous that he should have even dreamt of problems such as these, problems which it has taken the finest mathematicians in Europe a hundred years to solve, and of which the solution is incomplete to the present day.
\end{quote}
Furthermore, Watson \cite{gnwatson} wrote:
\begin{quote}
The most amazing thing about this formula is that it was discovered, apparently independently, by Ramanujan in his early days in India, and it appears in its appropriate place in his manuscript note-books.
\end{quote}

Among the 3200--3300 entries in his notebooks \cite{nb}, there are only a few instances in which Ramanujan returns to a topic that he had discussed elsewhere. In his third notebook \cite[Vol.~2, page 363, both editions]{nb}, Ramanujan reconsiders Claim \ref{ramaclaim}. Ramanujan's notebooks contain only a few indications of proofs, but more space is devoted to his   argument here than to any other argument or proof in his notebooks!  In particular, Ramanujan offers a \emph{more general theorem} than the one that  he employed to deduce Claim \ref{ramaclaim}.
In \cite[pp.~62--66]{Berndt} the first
author gave a detailed exposition of Ramanujan's argument
and, moreover, attempted to reconstruct his reasoning.

Put
$b(n)=1$ if $n$ can be written as a sum of two squares and $b(n)=0$, otherwise.
Since the time of Fermat, it has been well-known
that $b(n)=1$ precisely when in the canonical factorization of $n$, all primes
$p\equiv 3\,({\rm mod~}4)$ occur to an even exponent.  Thus, in particular, $b(n)$ is a
\emph{multiplicative} function.

Following Landau, we put $B(x)=\sum_{n\le x}b(n)$.
He proved in 1908
 that asymptotically
\begin{equation}
\label{Edmund}
B(x)\sim K\frac{x}{\sqrt{\log x}},\quad x\to\infty,
\end{equation}
where $K$ is given by (\ref{constantK}) (see \cite{L} or \cite[pp.\,641--669]{Lehrbuch}).
In honor of the contributions of both Landau and Ramanujan, the constant
$K$ is called the \emph{Landau--Ramanujan constant}. 
\pieter{In 1964, Shanks \cite{Shanks} determined 9 
decimals of $K$ using the fast converging expression
\begin{equation}
\label{Shanksconverging}    
K=\frac{1}{\sqrt{2}}\prod_{j=1}^{\infty}\Big(\frac{\zeta(2^j)(1-2^{-2^j})}{L(2^j)}\Big)^{(1/2)^{j+1}},\text{~with~}L(s)=\sum_{k=0}^{\infty}
\frac{(-1)^k}{(2k+1)^{s}},
\end{equation}
cf.\,\S \ref{sec:precise}.
Since 1964 the number of decimals of $K$ that can be computed
has increased dramatically.}
Indeed, D.~Hare \cite{Hare} claims to have computed
125,079  decimals; this result was recently superseded
by A.~Languasco \cite{Sconstant} who reached 130,000
 decimals.
In particular,
$$K=0.76422365358922066299069873125009232811679054\dotsc$$
For the latest on high precision evaluations of Euler
products involving primes in arithmetic progressions, such as the one for $K$, see S.~Ettahri et al.\,\cite{ERS} or 
O.~Ramar\'e \cite{Accurate}. A more leisurely account can be found in the book by Ramar\'e \cite[Chap.\,17]{walk}.


\par \pieter{Landau's proof of
\eqref{Edmund} uses complex analysis. In 1976, 
H.~Iwaniec \cite[Cor.\,1]{Iwaniec} gave a proof avoiding complex analysis based on sieving upper and lower bounds. A short and beautiful proof of 
Selberg \cite[pp.\,183--185]{Selberg} of an asymptotic result for the number of odd squarefree integers that can be written as a sum of two squares deserves honorable mention. Unfortunately, nobody was able to adapt this to a similar proof for $B(x)$.} 

\par Both Landau and Hardy were aware that
Landau's method can be extended
to show that $B(x)$ satisfies
an asymptotic
series expansion in the sense of Poincar\'e:
\begin{equation}
\label{serro}
B(x)=K\frac{x}{\sqrt{\log x}}\Bigl(1+\frac{c_0}{\log x}+\frac{c_1}{\log^2 x}
+ \dotsm + \frac{c_{m-1}}{\log^m x}
+O\Bigl(\frac{1}{\log^{m+1} x}\Bigr)\Bigr),\quad x\to\infty,
\end{equation}
where $m$ can be taken arbitrarily large and each $c_j$, $0\leq j \leq m-1$, is a constant.
Indeed, (\ref{serro}) became a folklore result, and a proof
was finally written down by J-P.~Serre \cite{Ser} for the larger class of so
called \emph{Frobenian multiplicative functions}.
These functions were subsequently considered
by R.W.K.~Odoni in a series of papers
(see \cite{O1} for a survey).
Serre \cite{Ser} also gave several beautiful applications  to coefficients of
modular forms, the origin of which goes back to Ramanujan's
Claim \ref{tauclaim}; see Section \ref{sec:claim2general}.
\par Put $$R(x)=K\int_2^x \frac{dt}{\sqrt{\log t}}.$$
Note that Ramanujan's claim implies, by partial integration
of $R(x)$, that
$$B(x)=K\frac{x}{\sqrt{\log x}}\Bigl(1+\frac{d_0}{\log x}+\frac{d_1}{\log^2 x}
+ \dotsm + \frac{d_{m-1}}{\log^m x}+O\Bigl(\frac{1}{\log^{m+1}x}
\Bigr)\Bigr),\quad x\to\infty,$$
 where $d_j=(2j+1)!/(j!\,2^{2j+1}).$
This seems promising, as asymptotically it is correct by (\ref{Edmund}) and the expansion follows the format (\ref{serro}).
However,
it turns out that $c_0\ne d_0=1/2,$ and thus Ramanujan's Claim \ref{ramaclaim} is
\emph{false}. This
was first asserted by
Gertrude Stanley \cite{Stanley}, a
Ph.D.\,student
of Hardy. Unfortunately, as
noted by D.~Shanks \cite{Shanks}, she made several errors (even after
her Corrigenda \cite{Stanley} are taken into account).
Stanley thought that $c_0<0$, and this led Hardy to the statement
\cite[p.\,19]{Hardy} that ``The integral is better replaced by the simpler
function $Kx/\sqrt{\log x}$''. This is false, as Shanks showed that
$c_0\approx 0.581948659$, which is very close to Ramanujan's predicted value $d_0=0.5$. 
Actually, Shanks \cite{Shanks} showed that Ramanujan's
integral numerically approximates
$B(x)$ closer than does $Kx/\sqrt{\log x}$ (for all $x$ large enough).
We will rederive his result (Theorem \ref{Shanksrecap}) 
as a particular case of something fitting in the more general framework discussed in Section \ref{sec:LvR}.
\par Shanks' constant $c_0$ can be computed these days with much higher precision: the
current record is by Languasco, who calculated it up to 130,000 decimal digits
\cite{Sconstant}. For a discussion of $c_0$ and related constants, see S.R.~Finch's book   \cite[Section 2.3]{constants}.

 It might be appropriate to end this introduction with what Hardy
wrote about Ramanujan in the context of Claim \ref{ramaclaim} being false \cite[p.~xxv]{cp}: ``And yet I am not sure that,
in some ways, his failure was not more wonderful than any of his
triumphs.''

\section{First Interlude: Gauss's Circle Problem}
\label{circleproblem}
Ramanujan's interest in sums of squares was not confined to Claim \ref{ramaclaim}.  Let $r_2(n)$ denote the number of representations of the positive integer $n$ as a sum of two squares, where different signs and different orders are counted as distinct.  For example, since
$$ 5=(\pm2)^2+(\pm1)^2=(\pm1)^2+(\pm2)^2,$$
$r_2(5)=8$. By identifying each representation of $n$ with a unit square, e.g.,  the square in which the lattice point lies in its southwest corner, it is easy to see that the number $R(x)$ of lattice points in a circle with radius $\sqrt{x}$ is approximately $\pi x$.  In other words,
\begin{equation}\label{PP}
 R(x)=\sum_{n\leq x}r_2(n) =\pi x+P(x),
 \end{equation}
where $r_2(0)=1$  and $P(x)$ is the error made in this approximation.  Observe that, from an examination of lattice points in circles of radii $\sqrt{x}+\sqrt2$ and $\sqrt{x}-\sqrt2$, we find that, respectively,
$$R(x)\leq (\sqrt{x} +\sqrt2)^2\qquad\text{and} \qquad
R(x)\geq (\sqrt{x} -\sqrt2)^2.$$
From these two inequalities, Gauss concluded that there exists a positive constant $C$ such that
\begin{equation}
\label{gaussp}|P(x)|\leq C\sqrt{x}\end{equation}
for all $x\geq0$.  Finding the optimal power $\theta$ such that $|P(x)|\leq C^{\prime}x^{\theta}$ for all $x>0$  is known as the \emph{Gauss Circle Problem}.
Hardy in his famous paper \cite{hardysquares} proved that
\begin{equation}\label{omega}
P(x)\neq O\bigl(x^{1/4}\bigr), \qquad  x\to \infty,
\end{equation}
and so in particular, $\theta>\frac14$.
In \cite{hardysquares}, Hardy offers a beautiful identity of Ramanujan that is not found elsewhere in Ramanujan's work, namely, if $a,b>0$, then
\begin{equation}\label{ab}
\sum_{n=0}^{\infty}\frac{r_2(n)}{\sqrt{n+a}}\,e^{-2\pi\sqrt{(n+a)b}} =
\sum_{n=0}^{\infty}\frac{r_2(n)}{\sqrt{n+b}}\,e^{-2\pi\sqrt{(n+b)a}}.
\end{equation}
Hardy used an identity that can be derived from (\ref{ab}) to prove (\ref{omega}).

Since Gauss derived \eqref{gaussp}, there have been many improvements on lowering the upper bound for $\theta$.  The most recent published record was set in 2003 by M.~Huxley \cite{huxley}, who proved that
$$\theta \leq \dfrac{131}{416}=0.3149\dotsc.$$
Recently, in an unpublished manuscript, X.~Li and Y.~Yang \cite{xx}, employing a very difficult argument, proved that $\theta\leq 0.31448\dotsc .$ 

The error term $P(x)$ can be represented by an infinite series of Bessel functions. More precisely,
\begin{equation}\label{circle}
{\sideset{}{^\prime} \sum_{0\leq n \leq x}}
r_2(n)= \pi x +
\sum_{n=1}^{\infty}r_2(n)\left(\frac{x}{n}\right)^{1/2}J_1(2\pi\sqrt{nx}),
\end{equation}
where the prime 
on the summation sign at the left indicates that if $x$ is an integer, then only $\frac12 r(x)$ is counted, and where $J_1(x)$ is the ordinary Bessel function of order 1.
The identity (\ref{circle}) apparently first appeared in Hardy's paper \cite{hardysquares},  where he wrote
 ``The form of
this equation was suggested to me by Mr.~S.~Ramanujan\dots.''  The identity \eqref{circle} is an important key in deriving upper bounds for $P(x)$.

In his Lost Notebook \cite[p.~335]{lnb}, Ramanujan offers a two-variable analogue of \eqref{circle}.  This and an analogous identity on page 335 of \cite{lnb}  involving $d(n)$, the number of positive divisors of the integer $n$, are the focus of several papers by the first author, S.~Kim, and A.~Zaharescu \cite{bkz1,bkz2,bkz3,bkz4,bkz5}.   These authors feel that Ramanujan derived the former identity to attack Gauss's \emph{Circle Problem}, but, although the first author, Kim, and Zaharescu have proofs of Ramanujan's two identities, they have been unable to penetrate Ramanujan's thinking.  For a discussion of these identities and an historical
overview of the \emph{Circle Problem}, see their survey article \cite{circledivisorsurvey1}.

\section{Second Interlude: Ramanujan's tau-function}
In non-technical terms, a
\emph{modular form} is a holomorphic function on the Poincar\'e upper half plane having a lot of symmetry with respect to the modular group $SL_2(\mathbb Z)$ or a congruence subgroup $\Gamma$ of it. In particular, for some fixed $k$, called the 
\emph{weight}, we have   
$$ 
f\Big(\frac{az+b}{cz+d}\Big)=(cz+d)^kf(z)
\quad
\text{~for~every~}
\left ( 
\begin{array}{cc}
 a  & b  \\
 c  & d  \\
\end{array}
 \right ) 
\in   \Gamma.
$$
We also require moderate growth at the cusps of the group of symmetry.
If the form vanishes at all of these cusps, it is said to be a \emph{cusp form}.
The lowest weight for cusp forms for the full modular group is 12, and it occurs for the \emph{normalized discriminant function $\Delta(z)$}. As already known to C.G.J.~Jacobi \cite{jacobi}, it can be expressed as an infinite product:
\begin{equation}
\label{simpleexpansion}
\Delta(z)=q\prod_{j=1}^{\infty}(1-q^j)^{24}=\sum_{n\ge 1}\tau(n)q^n,~~q=e^{2\pi iz},~~\Im(z)>0,
\end{equation}
called the modular discriminant. Its 24th root was extensively later studied by R.~Dedekind and is appropriately named the
\emph{Dedekind $\eta$-function}.
However, Ramanujan was the first to take a deep interest in arithmetic properties of $\tau(n)$ (now called \emph{Ramanujan's tau-function}),
and in 1916 published the important paper \cite{Rama1916}.
He calculated
$\tau(n)$ for $n=1,2,\dotsc,30$;
see Table \ref{tab:table1}. 
Note that the only odd values occur for $n=1,9,$ and
$25$ (see \S \ref{sec:parity} for a proof). Indeed, Ramanujan showed with ease that
$\tau(n)$ is odd
or even according as $n$ is an odd square or not.
He also made some observations that he was not able
to prove:

\begin{table}[ht]
\footnotesize
\renewcommand{\arraystretch}{1.25}
\begin{tabular}{|r|r||r|r||r|r|}
\hline
$n$ & $\tau(n)$ & $n$ & $\tau(n)$ & $n$ & $\tau(n)$ \\ \hline \hline
$1$ & $1$ & $11$ & $534612$ & $21$ & $-4219488$ \\ \hline
$2$ & $-24$ & $12$ & $-370944$ & $22$ & $-12830688$ \\ \hline
$3$ & $252$ & $13$ & $-577738$ & $23$ & $18643272$ \\ \hline
$4$ & $-1472$ & $14$ & $401856$ & $24$ & $21288960$ \\ \hline
$5$ & $4830$ & $15$ & $1217160$ & $25$ & $-25499225$ \\ \hline
$6$ & $-6048$ & $16$ & $987136$ & $26$ & $13865712$ \\ \hline
$7$ & $-16744$ & $17$ & $-6905934$ & $27$ & $-73279080$ \\ \hline
$8$ & $84480$ & $18$ & $2727432$ & $28$ & $24647168$ \\ \hline
$9$ & $-113643$ & $19$ & $10661420$ & $29$ & $128406630$ \\ \hline
$10$ & $-115920$ & $20$ & $-7109760$ & $30$ & $-29211840$ \\ \hline
\end{tabular}
\bigskip
\caption{The first 30 values of $\tau(n)$}
\label{tab:table1}
\end{table}
\begin{conjecture}
\label{fabulousproperties}
The following properties hold:\\
\rm{(a)}  $\tau(n)$ is multiplicative; that is, $\tau(mn)=\tau(m)\tau(n)$ whenever $(m,n)=1$;\\
\rm{(b)} if $p$ is prime, then $\tau(p^{e+1})=\tau(p)\tau(p^{e})-p^{11}\tau(p^{e-1})$
for any $e\ge2$;\\
\rm{(c)} $|\tau(n)|\le d(n)n^{11/2},$ where $d(n)$ denotes
the number of positive divisors of $n.$
\end{conjecture}

The first \pieter{and second property
imply} that the $\tau$-values are determined by their
values at prime arguments.
Ramanujan already noted that the first two properties can be rephrased in terms of the $L$-function 
$$L(\Delta,s) \coloneq \sum_{n\ge 1}\frac{\tau(n)}{n^s}=\prod_p \frac{1}{1-\tau(p)p^{-s}+p^{11-2s}},\quad
~~\Re(s)>\frac{13}{2}.$$
\pieter{Conjectures \ref{fabulousproperties}(a) 
and (b)
were} proved by L.J.~Mordell within the
year, and years later it motivated E.~Hecke \cite{Hecke} to introduce the \emph{Hecke operator} $T_p$ at each prime $p$ \pieter{leading to a rather more conceptual 
reproof of parts (a) and (b)}. In this set-up,  the product expansion is equivalent to $\Delta$ being an eigenfunction of $T_p$ with eigenvalue $\tau(p)$.
Using his operators, Hecke \cite{Hecke} went on to prove a similar result for Fourier coefficients of modular forms
for congruence subgroups of $\text{SL}_2(\mathbb Z)$.
\pieter{These days there is hardly any paper on modular forms that does not mention Hecke operators.}

Conjecture \ref{fabulousproperties}(c), however, resisted proof attempts by the best
minds for a long time.
For primes $p$,
 it states that in the factorization
$$1-\tau(p)p^{-s}+p^{11-2s}=\bigl(1-\alpha(p)p^{-s}\bigr)\bigl(1-\beta(p)p^{-s}\bigr),$$
we have
$|\alpha(p)|=|\beta(p)|=p^{11/2}$, and
hence
\begin{equation}
\label{Delignebound}
    |\tau(p)|\le |\alpha(p)|+|\beta(p)|\le 2p^{11/2}.
\end{equation}    
Ramanujan himself showed that
$\tau(n)=O(n^7)$ (we note that $d(n)=O(n^{\epsilon})$ for
any $\epsilon>0$), \pieter{later authors, 
cf.\,R.~Rankin \cite[\S 3]{rankintau}, used analytic number theory to improve on this, but did not come very close either in proving Conjecture \ref{fabulousproperties}(c).}
In 1971, P.~Deligne \cite{DelShi} interpreted the numbers $\tau(p)$ as eigenvalues of the Frobenius
automorphism on the cohomology of an appropriate 11-dimensional variety (generalization
of the Eichler-Shimura isomorphism) to reduce
the Ramanujan conjecture to the Weil conjectures for smooth projective varieties over finite fields,
which he
proved in 1974 \cite{Deligne} (see N.M.~Katz \cite{Katz} and E.~Kowalski \cite{Kowalski} for
introductory accounts). This proof is one of the most
important and highly regarded proofs in all
of 20th century algebraic geometry.

Conjecture \ref{fabulousproperties}(c) can be extended (in a suitable manner depending on the weight) to all so-called normalized Hecke eigencuspforms. R.P.~Langlands \cite{Langlands} reinterpreted the Hecke eigenforms in terms of automorphic representations for $GL_2$ over the rationals, so that both have the same associated $L$-functions. As such, the cusp forms correspond to cuspidal representations. The Ramanujan conjecture, in its general form, asserts that a generic cuspidal automorphic irreducible unitary representation of a reductive group over a global field should be locally tempered everywhere. It is largely unsettled, and research on it is part of the Langlands Program, a very active modern research area (see, S.~Gelbart \cite{Gelbart} for an introduction).
For more on the mathematical road from Conjecture \ref{fabulousproperties} to the Langlands Program, see the survey by W.-C.W.~Li \cite{winnie}
or the book by
the Murty brothers \cite{MM}.
The three excellent survey articles highlighting different aspects of the tau-function \cite{MurtyRR,rankintau,S-D} in
the proceedings of the 1987 ``Ramanujan Revisited'' conference, focus more on the properties of 
the tau-function per se.

It is amazing to see that the properties
Ramanujan uncovered and those that he conjectured  initially formed an easily overlooked trail,  but they are now  major mathematical highways!

\subsection{Congruences for \texorpdfstring{$\tau(n)$}{taun3}}
\label{taucongruences}
Ramanujan  also discovered congruences for $\tau(n)$, mostly involving the
\emph{sum of divisors} function $\sigma_k(n),$ which is defined as
$\sigma_k(n)=\sum_{d\mid n}d^k$.
Here we give a sampling (taken from \cite{S-D}), in which $(\frac{a}{b})$ denotes the \emph{Legendre symbol}:
\begin{equation}
\label{ramacongruences}
\tau(n)\equiv
\begin{cases}
\sigma_{11}(n) \pmod{2^8}, \ \text{if}\ 2\nmid n,\\
n^2\sigma_{7}(n) \pmod{3^3},\\
n\sigma_{9}(n)\pmod{5^2},\\
n\sigma_{3}(n)\pmod{7},\\
\sigma_{11}(n)\pmod{691},\\
0\pmod{23}, \ \text{if}\ (\frac{n}{23})=-1.\\
\end{cases}
\end{equation}

The \pieter{last of these congruences} was refined in 1930 by J.R.~Wilton \cite{Wilton}. If we restrict $n$ to
be a prime $p$, this refinement yields:
\begin{equation}
\label{Ramanujanwilton}    
\tau(p)\equiv
\begin{cases}
1\pmod{23},&\,\text{if}\ p=23,\cr
0\pmod{23},&\, \text{if}\ (\frac{p}{23})=-1,\cr
2\pmod{23},&\, \text{if}\ p=X^2+23Y^2\ \text{with}\ X\ne 0,\cr
-1\pmod{23},&\, \text{otherwise}.
\end{cases}
\end{equation}

The primes satisfying the \pieter{second congruence 
are primes in a finite union} of arithmetic progressions. For the
third congruence, this is no longer the case; it is inherently non-abelian. Wilton's starting point is 
the trivial observation that modulo 23  the $n$-th Fourier coefficient $t(n)$ of $\eta(z)\eta(23z)$ equals $\tau(n)$.

The congruence \eqref{Ramanujanwilton} was derived 
in a different way by F.~van der Blij \cite{Blij}, who considered the 
number of the representations of the integer $n$ by a form of class $F_i$, with $F_i$ a reduced form of discriminant $-23$ (of which there are three, namely $F_1=X^2+XY+Y^2$, $F_2=2X^2+XY+3Y^2$ and $F_3=2X^2-XY+3Y^2$). He 
establishes the formulae
$$a(n,F_1)=\frac{2}{3}\sum_{d\mid n}\Big(\frac{d}{23}\Big)+\frac{4}{3}t(n),\quad a(n,F_2)=a(n,F_3)=a(n,F_1)-2t(n),$$
and deduces the congruences from this.
By related arguments, cf.\,Serre \cite{SerreJordan},
it can be shown that
$$\tau(p)\equiv N_p(X^3-X-1)-1~({\rm mod~}23),$$ 
where
$N_p(f)$ denotes the number of distinct roots modulo $p$ of a polynomial $f\in \mathbb Z[x].$ The crux is that the $L$-function of $\eta(z)\eta(23z)$ is closely related to the Dedekind zeta-function of the cubic field of $\mathbb Q[X]/(X^3-X-1)$ of discriminant $-23$, the Galois closure of which is the Hilbert class field of $\mathbb Q(\sqrt{-23})$, cf.\,D.~Joyner \cite[pp.\,50-54]{Joyner}. A.~Akbari and Y.~Totani \cite{AT} recently generalized the result of van der Blij to the case of binary quadratic forms of discriminant $-D$ with odd $h(-D)$.

Another non-obvious result involves the \emph{Padovan sequence} \cite{OEIS},
which
is defined by $B_0=0, B_1=B_2=1$, and
by
$B_n=B_{n-2}+B_{n-3}$ for every $n\ge 3.$
The second author and A.~Noubissie \cite{MN} showed that
a prime $p$ divides
$B_{p-1}$ if and only if
$\tau(p)\equiv 2~({\rm mod~}23).$

After Wilton's paper many further ones with congruences for the 
tau-function appeared.
This whole rag-bag of
results begged for a more theoretical explanation,
which only became available in the
early 1970s, and rests
on the Serre--Deligne representation theorem and the theory
of modular forms modulo $\ell.$ Thus the congruence modulo 691 in \eqref{ramacongruences} is a consequence of
the $\Delta$-function reducing to the Eisenstein series $E_{12}$
modulo 691. 
The Eisenstein series have sums of divisor functions as Fourier coefficients.
For a detailed survey, we refer the reader to H.P.F.~Swinnerton-Dyer's paper \cite{S-D}. 
\pieter{In the modern approach to 
understanding these congruences mod $\ell$, the theory of Galois representations plays a major role, with interesting congruences typically arising if the image of the representation 
in $\text{GL}_2(\mathbb F_{\ell})$ 
is small. Focusing only on divisibility by $\ell$ of the Fourier coefficients allows us to rather consider the traceless matrices in the projectification $\text{PGL}_2(\mathbb F_{\ell})$. Here thanks to a celebrated result of Dickson all the possible images 
have been classified. For each of them 
certain information on the $\ell$-divisibility 
of the Fourier coefficient 
can be deduced, for details see, e.g., M.~Daas \cite{Daas}.}

\par Many of the results proved in this setting eventually
played a role in the work on Fermat's Last Theorem by Ribet, Wiles
and others! Thus here again
Ramanujan's trail would eventually become a highway.

\subsection{Computation of \texorpdfstring{$\tau(n)$}{taun2}-values}
As already mentioned, Ramanujan  computed the first 30 tau-values
in 1916 \cite{Rama1916} (see Table \ref{tab:table1}).
Watson \cite{watson3}, according to Rankin
\cite{rankinzelf}, as a pass-time during World War II, extended these computations to 1000 values.
Around the same time D.H.~Lehmer
\cite{Lehmertable} computed $\tau(n)$ for similar ranges.  The objectives
of these early tabulators were the
following: to check Conjecture \ref{fabulousproperties}(c), to find primes $p$ for which $p$ divides $\tau(p)$, and
to find $n$ for which $\tau(n)=0.$
Unsurprisingly, they never found a counter-example to
Conjecture \ref{fabulousproperties}(c).
It is now known that $p$ divides $\tau(p)$ for $p=2,3,5,7,2411$ and $p=7758337633,$ and that there are no
further $p<10^{10}$ with this property; see the paper by N.~Lygeros and O.~Rozier \cite{LR}. Using various congruences, Serre \cite{Serre73} proved that if $\tau(p)=0$, then $p$ lies in one of $33$ congruence classes modulo $3488033912832000=2^{14}\cdot 3^7\cdot 5^3\cdot 7^2\cdot 23\cdot 691$.
The conjecture that $\tau(n)\ne 0$ became famous as Lehmer's Conjecture. Any known
computational approach to this makes use of the known congruences to discard many integers $n.$ Lehmer \cite{Lehmer} himself got
to $3316798.$
The current record (see
\cite{DHZ}) is $816212624008487344127999\approx 8 \cdot 10^{23}.$ This
impressive result ultimately rests on sophisticated methods from algebraic geometry,
which allow one to compute $\tau(p)$  in polynomial time.
There is even an entire  book \cite{taucomp}
devoted to proving this result. The basic idea is to efficiently compute $\tau(p)$ modulo enough small primes $\ell$ such that the Chinese remainder theorem and the bound
$|\tau(p)|\le 2p^{11/2}$ allow
one to uniquely determine $\tau(p)$. For example, $19$ divides $\tau(10^{1000}+46227)$.
With the use of deep methods in analytic number theory, it was recently proved that the density of integers $n$ with $\tau(n)=0$ is at most $1.15\cdot 10^{-12}$ \cite{HIS}.

A strengthening of
Lehmer’s conjecture, was proposed by A.O.L.~Atkin and J-P.~Serre \cite{Ser}, who conjectured  that $|\tau(p)|\gg_{\epsilon} p^{9/2-\epsilon}$
for every $\epsilon > 0$.
This conjecture implies, in particular, that given any fixed integer $a$, there are at most finitely many primes $p$ for which
$\tau(p) = a$. If $a$ is odd, then Murty
et al.\,\cite{MMS} proved that
there are at most finitely many integers $n$ for which $\tau(n)=a$ (these $n$ must be odd
squares). More precisely, they demonstrated the existence of an effectively computable positive
constant $c$ such that if $\tau(n)$ is odd, then $|\tau(n)| > (\log n)^c$. M.A.~Bennett et al.\,\cite{oddies} proved results on the largest prime factor of odd $\tau(n)$ values.
J.S.~Balakrishnan et al.\,\cite{BCO} showed that
$\{\pm 1,\pm 3,\pm 5,\pm 7,\pm 691\}$ do not occur as $\tau$-values, and in the follow-up paper \cite{BCO2} many other small odd values are excluded. Currently there is a flurry of activity excluding further values. Meanwhile some further \emph{even} values can be excluded, e.g., $2p$, where $p$ is a prime,
$2<p<100$ or $-2p^j$ with $j\ge 1$ arbitrary and
with $p$ in a set containing eighteen primes \cite{BOT}.

Serre \cite{SerreChebotarev} initiated the general study of estimating the size of possible gaps in the Fourier expansions of modular forms, namely, he considered the gap function
$$i_f(n)=\max\{k:a_f(n+j)=0\text{~for~all~}0\le j\le k\}.$$
For Fourier coefficients of a newform $f$ without complex multiplication, A.~Balog and K.~Ono \cite{BalogOno} showed that
$$i_f(n)\ll_{f,\epsilon}n^{\frac{17}{41}+\epsilon}.$$
E.~Alkan \cite{Alkan}
improved the exponent $17/41$ to $51/134$ in case of newforms associated to elliptic curves without complex multiplication. 
In \cite{Alkan2} he showed that for these newforms $i_f(n)^{\lambda}$ is bounded on average for $\lambda<1/8$.
Alkan and A.~Zaharescu \cite{AZ}, making clever use of the
Ramanujan--Wilton congruence \eqref{Ramanujanwilton}, showed that $i_{\Delta}(n)\le 2\sqrt{46}n^{1/4}+23$ for every $n\ge 1$. By a variation of their method S.~Das and S.~Ganguly \cite{DG} showed that $i_{f}(n)\ll_k n^{1/4}$ for any non-zero cusp form $f$ of integral 
weight $k$ on the full modular group.
Interestingly, in the final step of their analysis they have to
establish the existence of sums of two squares in short intervals.

\section{Third Interlude: Euler--Kronecker constants of multiplicative sets}
\label{sec:LvR}
A set $S$ of natural numbers is said to be {\it multiplicative} if for every pair $m$ and $n$ of co-prime integers in $S$,
$mn$ is also   in $S$.
Let $i_S$ denote the
\emph{characteristic function} of $S$, i.e.,
\begin{equation*}
i_S(n)=\begin{cases} 1, \quad & n\in S,\\
0, &\textup{otherwise}.
\end{cases}
\end{equation*}
 Note that the set $S$ is multiplicative if and only if $i_S$ is a multiplicative function.
A large class of multiplicative sets is provided by the sets
\begin{equation}
    \label{setdefinition}
S_{f;q} \coloneq \{n:~q\nmid f(n)\},
\end{equation}
where $q$ is any prime and $f$ is any integer-valued multiplicative function.

Let $\pi_S(x)$ and $S(x)$ denote the number of primes and, respectively, the number of integers in $S$ not exceeding $x.$
The following result is a special case of a theorem
of E.~Wirsing \cite{Wirsing}, with a reformulation following Finch et al.\,\cite[p.\,2732]{FMS}. As usual,
$\Gamma$  denotes the gamma function.

\begin{theorem}
\label{een}
Let $S$ be a multiplicative set satisfying $\pi_S(x)\sim \delta x/\log x$, as $x\to\infty$,  for some $\delta$, $0<\delta<1$.
Then $$S(x)\sim c_0(S) x \log^{\delta-1}x,\quad x\to\infty, $$
where
$$c_0(S) \coloneq \frac{1}{\Gamma(\delta)}\lim_{P\rightarrow \infty}\prod_{p<P}
\Bigl(1+\frac{i_S(p)}{p}+\frac{i_S(p^2)}{p^2}
+\dotsm \Bigr)\Bigl(1-\frac{1}{p}\Bigr)^{\delta}$$
converges and hence is positive.
\end{theorem}

Recall that for any character $\chi$,  Dirichlet's $L$-function $L(s,\chi)$ is defined by
$$L(s,\chi) \coloneq \sum_{n=1}^{\infty}\dfrac{\chi(n)}{n^{s}}, \quad \Re(s)>1.$$

\begin{example}
If $S$ is the set of integers that can be written as a sum of two squares, Theorem \ref{een} yields
$$K=c_0(S)=\frac{\sqrt{2}}{\Gamma(1/2)}\Bigl(L(1,\chi_{-4})\prod_{p\equiv3\pmod{4}}
\frac{1}{1-p^{-2}}\Bigr)^{1/2},$$
where $\chi_{-4}$  denotes the nontrivial, quadratic Dirichlet character
modulo $4$, and where we used the identity
$$\lim_{P\rightarrow \infty}\prod_{\substack{p<P\\ p\equiv 1\pmod{4}}}
\Bigl(1-\frac{1}{p}\Bigr)^{-1}\prod_{\substack{p<P\\ p\equiv 3\pmod{4}}}
\Bigl(1+\frac{1}{p}\Bigr)^{-1}=L(1,\chi_{-4}).$$
Since $\Gamma(1/2)=\sqrt{\pi}$ and
$L(1,\chi_{-4})=\pi/4$, we then obtain the expression \eqref{constantK} for $K$.
\end{example}

In this set-up, Ramanujan likely would have claimed that
\begin{equation}
\label{valseanalogiegeneral}
S(x)=c_0(S)\int_2^x \log^{\delta-1}t\, dt +O(x\,\log^{-r} x),
\end{equation}
where $r$ is any positive number.
We call $$c_0(S) x \log^{\delta-1}x,~~c_0(S)\int_2^x  \log^{\delta-1}t\,dt,$$
respectively, the \emph{Landau} and
\emph{Ramanujan approximations}
to $S(x)$.
If for all $x$ sufficiently large, $$\Big|\,S(x)- c_0(S) x \log^{\delta-1}x\,\Big|<
\Big|\,S(x)- c_0(S) \int_2^x  \log^{\delta-1}t\,dt\,\Big|,$$
 we say that the Landau approximation is better than the
Ramanujan approximation. If the reverse inequality holds for every $x$ sufficiently large, we say
that the Ramanujan approximation is better than the Landau approximation.
We will now introduce  a constant which can be used to decide which approximation is better.

For $\Re(s)>1$, put
$$L_S(s) \coloneq \sum_{n\in S}n^{-s}.$$
If the limit
\begin{equation}
\label{EKf}
\gamma_S \coloneq \lim_{s\rightarrow 1^+}\Bigl(
\frac{L'_S(s)}{L_S(s)}
+\frac{\alpha}{s-1}\Bigr)
\end{equation}
exists for some $\alpha>0$, we say that the set $S$ admits an \emph{Euler--Kronecker constant} $\gamma_S$. In the case
$S=\mathbb N,$ we have
$L_S(s)=\zeta(s)$, the
\emph{Riemann zeta function}, $\alpha=1$ and $\gamma_S=\gamma=0.5772156649\dotsc$, the \emph{Euler--Mascheroni
constant}. Consult J.~Lagarias's paper \cite{Lagarias} for a beautiful  survey, and
G.~Havil \cite{Havil} for a popular account. If $S$ is a set that in some sense is close to the set of all natural numbers, then $\gamma_S$ will be close to
$\gamma$. This will be, for example,
the case if $q$ in \eqref{setdefinition} is a large prime (see, e.g., Example \ref{Fordexample}).

As
the following result shows, the Euler--Kronecker constant $\gamma_S$ determines the second order behavior of $S(x).$ As usual $\pi(x)$ denotes the prime counting function.

\begin{theorem}
\label{thm:multiplicativeset}
Let $S$ be a multiplicative set.
If there exists $\rho>0$ and $0<\delta<1$ such that
\begin{equation}
\label{primecondition}
\pi_S(x)=
\delta\pi(x)+O_S(x\log^{-2-\rho}x),\quad x\to\infty,
\end{equation}
then $\gamma_S\in\mathbb R$ exists, and asymptotically, as $x\to\infty$,
\begin{equation}
\label{initstarrie}
S(x)=\sum_{\substack{n\le x \\n\in S}}1=\frac{c_0(S)\,x}{\log^{1-\delta}x}\left(1+\frac{(1-\gamma_S)(1-\delta)}{\log x}\bigl(1+o_S(1)\bigr)\right).
\end{equation}
In the case when  the prime numbers belonging to $S$ are, with finitely many exceptions, precisely
those in a finite union of arithmetic progressions,
we have, for arbitrary $j\ge 1$,
\begin{equation}
\label{starrie1}
S(x)=\frac{c_0(S)\,x}{\log^{1-\delta}x}\Bigl(1+\frac{c_1(S)}{\log x}+\frac{c_2(S)}{\log^2 x}+\dotsm+
\frac{c_j(S)}{\log^j x}+O_{j,S}\Bigl(\frac{1}{\log^{j+1}x}\Bigr)\Bigr),\quad x\to\infty,
\end{equation}
where  $c_1(S)=(1-\gamma_S)(1-\delta)$, and
$c_2(S),\dotsc,c_j(S)$ are further constants.
\end{theorem}

\begin{proof}
\pieter{The first assertion is a consequence of Theorem 1.2 of
de la Bret\`eche and Tenenbaum
\cite{BT3},} see also Moree \cite[Theorem 4]{Mpreprint}; for the second, see Serre \cite[Th\'eor\`eme 2.8]{Ser}.
\end{proof}

By partial integration, as $x\to \infty$,
\begin{equation}
\label{part1}
\int_2^x  \log^{\delta-1}t\,dt=\frac{x}
{\log^{1-\delta}x}\Bigl(1+\frac{1-\delta}{\log x}+O\Bigl(\frac{1}{\log^2 x}\Bigr)\Bigr).
\end{equation}
Thus, the Landau and Ramanujan approximations to $S(x)$ amount to taking $c_1(S)=0$ and $c_1(S)=1-\delta$, respectively.
This trivial observation leads to the following corollary of Theorem \ref{thm:multiplicativeset}.

\begin{corollary}
\label{gammacomparison}
Suppose that the set $S$ is multiplicative and satisfies \eqref{primecondition}. Then
the Euler--Kronecker constant $\gamma_S$ exists.  Furthermore,
\begin{itemize}
\item The Ramanujan type claim \eqref{valseanalogiegeneral} is true for every $r\le 2-\delta$, and, provided that $\gamma_S\ne 0$, is false for every $r>2-\delta$.
\item If $\gamma_S>1/2$, the Landau approximation is better. If $\gamma_S<1/2$, the Ramanujan approximation is better.
\end{itemize}
\end{corollary}

\begin{example}
\label{Fordexample}
K.~Ford et al.\,\cite{FLM} studied the infinite family of sets
$S_{\varphi;q} \coloneq \{n:~q\nmid \varphi(n)\},$ where $q\ge 3$ is a prime and $\varphi$ denotes Euler's totient function. They showed that the Ramanujan approximation is better if and only if $q\le 67$. Further, they established that, as $q\to\infty$,
$\gamma_{S_{\varphi;q}}=
\gamma+O_{\epsilon}(q^{\epsilon-1})$, underlining the fact that for large $q$, the series $L_{S_{\varphi;q}}(s)$ starts to behave more like $\zeta(s)$. They also show that 
$\gamma_{S_{\varphi;q}}$ is
intimately related to the Euler--Kronecker constant of the cyclotomic number field $\mathbb Q(\zeta_q)$. In general the Euler--Kronecker constant of a number field $K$ is obtained on putting $\alpha=1$ and $L_S(s)=\zeta_K(s)$ in \eqref{EKf}, 
where $\zeta_K(s)$ denotes the Dedekind zeta function of $K$. 
\end{example}

\subsection{Abelian multiplicative subsets}
\label{abeliansubset}
The sets in Example \ref{Fordexample} and many others in the literature under the umbrella of Theorem \ref{thm:multiplicativeset} are abelian, which we now define.
\begin{definition}
A multiplicative set $S$ is called \textbf{abelian} if 
it consists, with
finitely many exceptions, of all the primes in a finite union of arithmetic progressions.
\end{definition}
\noindent For example, the set consisting of $5,7$ and all primes $p$ satisfying
$p\equiv \pm 1\pmod{7}$ is abelian. 

By the Chinese remainder theorem, we can find an integer $d$ such that the primes are, with finitely many exceptions, those in a number of primitive residue classes modulo $d$.

If $S_1$ and $S_2$ are multiplicative sets such that any two elements of $S_1$ and $S_2$ are co-prime, then
$S_1\cdot S_2 \coloneq \{m\cdot n:m\in S_1,\,n\in S_2\}$ is a multiplicative set and
$L_{S_1\cdot S_2}(s)=L_{S_1}(s)L_{S_2}(s)$.
If both Euler--Kronecker constants $\gamma_{S_1}$ and $\gamma_{S_2}$ exist, then $\gamma_{S_1\cdot S_2}=\gamma_{S_1}+\gamma_{S_1}$. Thus, in the case when $S$ is multiplicatively abelian, the computation of $\gamma_S$ can be reduced to the case where the primes in $S$ are precisely those in \emph{one} primitive residue class.

{}From now on we use the short-hand notation defined on the right-hand side below:
$$\frac{d}{ds}\log\{A(s)\}=\frac{A^{\prime}(s)}{A(s)} \eqcolon \frac{A^{\prime}}{A}(s).$$

\begin{theorem}[Languasco and Moree \cite{LaMoree}]
\label{mainabelian}
Let $a$ and $d\ge 2$ be co-prime integers. Let $S$ be the set of integers including $1$ and all integers composed only of primes $p\equiv a\pmod{d}$. Then $S$ has an Euler--Kronecker constant $\gamma(d,a)$ given by
\begin{equation}
    \label{compact}
\gamma(d,a)=\gamma_1(d,a)-\sum_{p\equiv a\pmod{d}}\frac{\log p}{p(p-1)}
+
\sum_{n\equiv a\pmod{d}}\frac{(1+\mu(n))\Lambda(n)}{n},
\end{equation}
where $\mu$ is the M\"obius function,
$\Lambda$ is the von Mangoldt function,
\begin{equation}
\label{gamma1-def}
\gamma_1(d,a)
=
\frac{1}{\varphi(d)}\Bigl(\gamma+\sum_{p\mid d}\frac{\log p}{p-1}+
\sum_{\chi\ne \chi_0}{\overline{\chi}}(a)
\frac{L'}{L}(1,\chi)\Bigr),
\end{equation}
and $\chi_0$ is the principal character modulo $d$.
\end{theorem}
It follows from this result that an abelian multiplicative set has an Euler--Kronecker constant involving Dirichlet $L$-series.

\begin{remark}
\label{Landau1909}
In 1909, using the method with which he proved the asymptotic formula \eqref{Edmund}, Landau \cite{La09a}
(see also \cite[\S 176--183]{Lehrbuch})
 settled a question
of D.N.~Lehmer who, reformulated in our terminology, asked
about the asymptotic behavior of
$$x^{-1}\sum_{\substack{n\le x,\\n\in S}}2^{\omega(n)},$$
where $\omega(n)$ denotes the number of different prime factors of $n$, and $S$ is assumed to be abelian. In this context, he established \eqref{starrie1}, however, without identifying
$c_1(S)$ as $(1-\gamma_S)(1-\delta)$.
\end{remark}

\subsection{The Euler--Kronecker constant for sums of two squares}
\label{sec:EKsumoftwosquares}
As an example, we will determine $\gamma_S$ in the case when $S$ is the set of sums of two squares.
Since $S$ is generated by the prime $2$, the primes $\equiv 1\pmod{4}$ and the squares of the primes $\equiv 3\pmod{4}$, we obtain 
\begin{equation}
 \label{LSdefine}   
L_S(s)=(1-2^{-s})^{-1}\prod\limits_{p\equiv 1 \pmod{4}}(1-p^{-s})^{-1}
\prod\limits_{p\equiv 3 \pmod{4}}(1-p^{-2s})^{-1}.
\end{equation}
Recall that
$$L(s,\chi_{-4})=\prod\limits_{p\equiv 1 \pmod{4}}(1-p^{-s})^{-1}
\prod\limits_{p\equiv 3 \pmod{4}}(1+p^{-s})^{-1}.$$
Comparing the Euler factors on both sides, we can verify that
\begin{equation}
    \label{LSsquare}
L_S(s)^2=\zeta(s)L(s,\chi_{-4})(1-2^{-s})^{-1}
\prod\limits_{p\equiv 3 \pmod{4}}(1-p^{-2s})^{-1}.
\end{equation}
For the reader who is familiar with the \emph{Dedekind zeta function} $\zeta_K(s)$, this identity is not so mysterious, and will realize that
$$\zeta_{\mathbb Q(i)}(s)=\zeta(s)L(s,\chi_{-4})=(1-2^{-s})^{-1}\prod\limits_{p\equiv 1 \pmod{4}}(1-p^{-s})^{-2}
\prod\limits_{p\equiv 3 \pmod{4}}(1-p^{-2s})^{-1},$$
which exhibits the relationship of $L_S(s)^2$ with $\zeta_{\mathbb Q(i)}(s)$.

Applying logarithmic differentiation to both sides of \eqref{LSsquare}, we find that 
$$2\frac{L_S'}{L_S}(s)+\frac{1}{s-1}=\frac{\zeta'}{\zeta}(s)+\frac{1}{s-1}+\frac{L'}{L}(s,\chi_{-4})-\frac{\log 2}{2^s-1}-\sum\limits_{p\equiv 3 \pmod{4}}\frac{2\log p}{p^{2s}-1}.$$
Therefore, letting $s\to1$, we deduce that
\begin{equation}
\label{theanswer}
2\gamma_S=\gamma+\frac{L'}{L}(1,\chi_{-4})-\log 2-\sum_{p\equiv 3\pmod{4}}\frac{2\log p}{p^{2}-1}.
\end{equation}

\subsubsection{Rough numerical evaluation of $\gamma_S$}
The quotient
\begin{equation}\label{LL}
 \frac{L'}{L}(1,\chi_{-4})
 \end{equation}
in \eqref{theanswer} is connected with two important ideas in the mathematical literature: 
the arithmetic-geometric mean (AGM) of Gauss and Lagrange, and the lemniscate integral.

First, we discuss the AGM.  Let $a=a_0$ and $b=b_0$ be initial values, with $a_0, b_0>0$.  Recursively define two 
sequences \pieter{$(a_n)$ and $(b_n)$} for $n\geq 1$ by
\begin{equation*}
a_{n+1}=\dfrac{a_n+b_n}{2} \quad \textup{and} \quad b_{n+1}=\sqrt{a_nb_n}.
\end{equation*}
Then
\begin{equation*}
\lim_{n\to\infty}a_n  \quad \textup{and} \quad \lim_{n\to\infty}b_n
\end{equation*}
both exist and are equal.  The \emph{arithmetic-geometric mean} of 
\pieter{$(a_n)$ and $(b_n)$} is defined by
\begin{equation}\label{M}
M(a,b) \coloneq \lim_{n\to\infty}a_n =  \lim_{n\to\infty}b_n.
\end{equation}
To see how extensively the AGM appears in number theory and related analysis, consult J.M.~and P.B.~Borwein's fascinating treatise \cite{BoBo}.

Second, the \emph{lemniscate integral}, naturally arising in the calculation of the arc length of the lemniscate, is defined by
\begin{equation}\label{L}
L \coloneq \int_0^1\df{dx}{\sqrt{1-x^4}}.
\end{equation}
The lemniscate integral was initially studied by Count Giulio Fagnano and James Bernoulli.  C.L.~Siegel \cite{siegel} considered the lemniscate integral so important that he began his series of lectures on elliptic functions with a thorough discussion of it. For an interesting historical account of the lemniscate integral, including the work of Fagnano and Bernoulli, consult R.~Ayoub's paper \cite{ayoub}; see also a paper by G.~Almkvist and the first author \cite{AB}.

On pages 283, 285, and 286 in the unorganized portion of his second notebook \cite{nb}, Ramanujan examined the lemniscate integral  \eqref{L} and various extensions and analogues of it. In particular, he established several inversion formulas.  We state one of his results, which was first proved by S.~Bhargava and the first author \cite{bhargava1}.

\begin{entry} Let $\theta, v$ and $\mu$ be defined by
\begin{equation}\label{elliptic}
\df{\theta\mu}{\sqrt{2}}=\int_0^v\df{dx}{\sqrt{1-x^4}},
\end{equation}
where $0\leq \theta\leq \pi/2, 0\leq v \leq 1,$ and $\mu$ is a constant defined by setting $\theta=\pi/2$ and $v=1$.  Then, for $0< \theta\leq \pi/2$,
\begin{equation}\label{inversion}
\df{\mu^2}{2v^2}=\csc^2\theta-\df{1}{\pi}-8\sum_{n=1}^{\infty}\df{n\,\cos(2n\theta)}{e^{2\pi n}-1}.
\end{equation}
\end{entry}
We will deduce the value of the constant $\mu$ below.
For proofs of Ramanujan's inversion formulas, see the first author's book \cite[Chapter 26]{Berndt} or \cite{bhargava1}.

Returning to \eqref{theanswer} and \eqref{LL}, we 
\pieter{point out}
the identity
\begin{equation}
\label{classicalGauss}
\frac{L'}{L}(1,\chi_{-4})=\log\bigl(M(1,\sqrt{2})^2e^{\gamma}/2\bigr),
\end{equation}
where $M(1,\sqrt{2})$ is given by \eqref{M}  with starting values $a_0=1$ and $b_0=\sqrt{2}$. 
This \pieter{identity} appears to have been discovered independently
at least by Berger (1883), Lerch (1897), de S\'eguier (1899) and Landau (1903) (for more complete references, see Shanks \cite{Shanks}).
Gauss showed in his diary \cite{Gray} that
\begin{equation}\label{g}
G \coloneq \frac{1}{ M(1,\sqrt{2})}=\frac{2}{\pi}\int_0^1 \frac{dx}{\sqrt{1-x^4}}=0.8346268416740731862814297\dotsc.
\end{equation}
This constant is now named
\emph{Gauss's constant} and is
also related to various other constants. Note that when we compare \eqref{elliptic} with \eqref{g}, we deduce that \cite{Gaussconstant}
$$\mu=\df{\sqrt2}{M(1,\sqrt2)}.$$

We conclude from \eqref{classicalGauss} and \eqref{g} that
\eqref{theanswer} can be rewritten as
\begin{equation}
\label{theanswer2}
\gamma_S=\gamma-\log G-\log 2-\sum_{p\equiv 3\pmod{4}}\frac{\log p}{p^{2}-1}.
\end{equation}
The AGM algorithim is fast. Just taking a few decimals of the constants involved into account and only the prime $q=3$ in the sum, we conclude that
$\gamma_S<0.578-\log(5/6)-0.693-(\log 3)/8<-0.07$, where we note that it takes only two steps in the AGM algorithm in order to conclude that $G>5/6$. On applying Corollary
\ref{gammacomparison} we obtain a new proof of the following result.

 \begin{theorem}[Shanks \cite{Shanks}]
 \label{Shanksrecap}
Ramanujan's approximation $K\int_2^x \frac{dt}{\sqrt{\log t}}$
asymptotically better approximates $B(x)$, than does Landau's asymptotic $K \frac{x}{\sqrt{\log x}}$.
However, Ramanujan's Claim \ref{ramaclaim} is
false with any error term 
$\theta(x)$ satisfying $\theta(x)=o(x\log^{-3/2}x)$.
\end{theorem}

\subsubsection{Precise numerical evaluation of $\gamma_S$}
\label{sec:precise}
We first describe the computation of Shanks \cite{Shanks} \pieter{of both $K$ and $\gamma_S$}.
His starting point is formula \eqref{theanswer} and the observation that, for $\Re(s)>1/2$,
\begin{equation}\label{above}
\prod\limits_{p\equiv 3 \pmod{4}} (1-p^{-2s})^{-2}=\zeta(2s)(1-2^{-2s})L(2s,\chi_{-4})^{-1}\prod\limits_{p\equiv 3 \pmod{4}}(1-p^{-4s})^{-2}.
\end{equation}
\pieter{By induction one can infer from 
this \eqref{Shanksconverging}, with 
of course $L(s)=L(s,\chi_{-4})$.} 

Using the logarithmic differentiated form of 
\eqref{above} repeatedly, he obtains
\begin{equation}
\label{idie}
\sum\limits_{p\equiv 3 \pmod{4}}
\frac{2\log p} {p^2-1}=\sum_{k=1}^{\infty}\Bigl(\frac{L'}{L}(2^k,\chi_{-4})-\frac{\zeta'}{\zeta}(2^k)-\frac{\log 2}{2^{2^k}-1}\Bigr),
\end{equation}
which in combination with \eqref{theanswer}
and \eqref{classicalGauss} then yields
$\gamma_{S}=-0.1638973186\dotsc$.

These days, one can do much better.  To compute
$\frac{L'}{L}(j,\chi)$ for $j\ge 2$, first recall that $L(s,\chi)$ can be expressed as a linear combination of Hurwitz zeta functions $\zeta(s,x)$.
An efficient algorithm to compute $\zeta(s)$, $\zeta(s,x)$, and $\frac{d\zeta}{ds}(s,x)$ for $s>1$, has been devised by
 Languasco
 \cite{Languasco2023a}.
Invoking \eqref{theanswer} and \eqref{idie} he determined 
130,000 decimals for $\gamma_S$ in \cite{Sconstant}.

We note that the typical prime sums can be naively estimated
by computing all terms up to some large value $x$ and estimating the tail by
$$\sum_{p>x}\frac{\log p}{p^k-1}\le \frac{x}{x^k-1}\Bigl(-0.98+1.017\frac{k}{k-1}\Bigr),\quad k\in \mathbb R_{>1}\text{~and~}x\ge 7481,$$
which follows easily on using the estimate $0.98x \le \sum_{p\le x}\log p \le 
1.017x$ for $x \ge 7481$ due to J.B.~Rosser and L.~Schoenfeld \cite{RS}. The same bound can be used if the primes $p$ are restricted to some arithmetic progression.

\subsection{Connection with Cilleruelo's constant}
\label{sec:Cilleruelo}
The result mentioned in the previous section allows us to compute \emph{Cilleruelo's constant}
$$J \coloneq \gamma-1-\frac{\log 2}{2}-\sum_{p>2}\frac{(\frac{-1}{p})\log p}{p-1}=-0.0662756342\dotsc,$$
with many more decimals.\footnote{In the first formula given in Moree \cite{LvR} for $J$ read $p>2$ instead of $p>3$.}

Given an integer valued polynomial $f(X)$, we let
$L(N)$ denote the least common multiple of $f(1),\dotsc,f(N)$. In case $f(X)=X$ an equivalent
formulation of the Prime Number Theorem states that $\log L(N)\sim N$.
It is a natural question to generalize this to other polynomials.
For example, in case $f$ is a product of linear terms, 
we have $\log L(N)\sim c_f N$ for some positive constant $c_f$; see S.~Hong 
et al.\,\cite{HQT}. The
main ingredient of the proof is Dirichlet’s theorem on primes in arithmetic
progressions.
J.~Cilleruelo \cite{C} considers the case where $f$ is an irreducible polynomial
of degree 2 and shows that 
$\log L(N ) \sim N \log N$ and
conjectures that if $f$ is an irreducible polynomial of degree $d > 2$,
then 
$$\log L(N) \sim  (d - 1)N \log N.$$ Z.~Rudnick and S.~Zehavi \cite{RZ} showed that the conjecture is true on average for $f(X)-a$, with $a$ taken in a sufficiently large range.
J.~Maynard and Z.~Rudnick \cite{MR} 
proved that $$(1/d+o(1))N \log N\le \log L(N)\le (d-1+o(1))N\log N.$$
A.~Sah \cite{Sah} later established the  
lower bound with $1/d$ replaced by 1.

In case $f(X)=X^2+1$, Cilleruelo established the more precise result
\begin{equation}
\label{Javier}
\log L(N)=N\log N+JN+o(N).
\end{equation}
It is easy to see  that (\cite[Prop.\,4]{LvR})
$$-\sum_{p>2}\frac{(\frac{-1}{p})\log p}{p-1}=\frac{L'}{L}(1,\chi_{-4})+\sum_{p\equiv 3\pmod{4}}\frac{2\log p}{p^{2}-1}.$$
This, in combination with the identities \eqref{classicalGauss} and \eqref{theanswer2}, then yields
$$J=4\gamma-1-(7/2)\log 2-4\log G-2\gamma_S,$$
revealing
that from a computational point of view there is almost no difference between $J$ and $\gamma_S$, 
and also 130,000 decimals precision for $J$ can be obtained. It is an open problem to determine whether or not there exists a multiplicative set $C$ such that $\gamma_C=J$. 
It is believed, see Cilleruelo et al.\,\cite[p.\,104]{CRSZ},  that a result of the form \eqref{Javier} also holds for other irreducible quadratic monic polynomials with the analogue of $J$ being non-zero.

\section{Ramanujan and the non-divisibility of \texorpdfstring{$\tau(n)$}{taun}}
\label{nondivtau}
Let us return to Ramanujan's unpublished, highly influential manuscript on $\tau(n)$ and the partition function $p(n)$, the number of ways of representing a positive integer $n$ as a sum of positive integers, irrespective of their order. Portions of the manuscript were likely written in the years 1917--1919, prior to Ramanujan's return to India in 1919, while other parts might have been written after he returned to India.  Most likely, the manuscript was sent to Hardy by Francis Dewsbury, Registrar at the University of Madras, in 1923. This shipment of papers contained several unpublished manuscripts and fragments of Ramanujan, including what was later to be  called, \emph{Ramanujan's Lost Notebook} \cite{lnb}.  The manuscript is in two parts.  The first is 43 pages long and is in Ramanujan's handwriting; the second is 6 pages long and is in the handwriting of Watson.  Unfortunately, the second portion in Ramanujan's handwriting has never been found.  (There are several manuscripts of Ramanujan that exist only in Watson's handwriting.  Evidently, he copied them for his own use and then, sadly, discarded Ramanujan's original manuscripts.)  Supplying details when needed, the first author and Ono \cite{BerndtOno} made a thorough examination of the manuscript. A revised version of their study appears in \cite[Chapter 5]{lnbIII}. Congruences for $\tau(n)$ and $p(n)$ are highlights of the manuscript.

 In this manuscript, Ramanujan considers, for the primes $q=3,5,7,23,691,$
appearing in his congruences \eqref{ramacongruences}, the quantity
\begin{equation}
\sum_{n\le x,\, q\nmid \tau(n)}1
\end{equation}
and makes claims similar to Claim \ref{ramaclaim}.
He defines
\begin{equation*}
t_n =\begin{cases}
1,  \quad&\text{if } q\nmid \tau(n),\\
0,  &\text{otherwise},
\end{cases}
\end{equation*}
and then \pieter{(for some $q\in \{3,5,7,23,691\}$)}  writes:
\begin{claim}
\label{tauclaim}
It is easy to prove by quite
elementary methods that $\sum_{k=1}^n t_k=o(n)$, as $n\to\infty$. It can be shown by
transcendental methods that
\begin{equation}
\label{ditgaatnoggoed}
\sum_{k=1}^n t_k\sim \frac{C_qn}{(\log n)^{\delta_q}};
\end{equation}
and
\begin{equation}
\label{valseanalogie}
\sum_{k=1}^n t_k=C_q\int_1^n \frac{dx}{(\log x)^{\delta_q}}+O\Bigl(\frac{n}{(\log n)^r}
\Bigr),
\end{equation}
where $r$ is any positive number.
\end{claim}
Ramanujan used the notation $C$ and $\delta,$ but for us, to indicate their dependence on $q$, it is
more convenient to use $C_q$ and $\delta_q.$
\pieter{The} values of $\delta_q$ are given in the final
column of Table \ref{tab:table3}.
Note that the truth of $\sum_{k=1}^n t_k=o(n)$
would imply that $q|\tau(n)$ for almost all $n$.

 It appears from Stanley's paper \cite{Stanley}, which we discussed earlier in Section \ref{Ramasquare}, that Hardy planned to have this manuscript published under Ramanujan's name after some editing. Although he
published some parts of it (see \cite{BerndtOno}), unfortunately, he never published Ramanujan's
full manuscript. Proofs of some of Ramanujan's  assertions in his unpublished
manuscript were further worked out by  Stanley \cite{Stanley}. She asserted \eqref{valseanalogie} to be false for
$q=5$ and any $r>1+\delta_5$ with $\delta_5=1/4$
(cf.\,Table \ref{tab:table3}).
Corrections to her paper have been given by
the second author \cite[Section 5]{MRama}.

 In 1928, Hardy passed on
the unpublished manuscript to Watson, who unfortunately kept it hidden away from
the mathematical community. (F.J.~Dyson colorfully described Watson's penchant for keeping things to himself \cite{Dyson}.)
Watson wrote approximately 30 papers devoted to Ramanujan's work; see Rankin \cite{watsonobi} for an overview.

Note that
Ramanujan's Claims \ref{ramaclaim} and \ref{tauclaim}
 are very reminiscent of the Prime Number Theorem in the form
\begin{equation}
\label{hadamard}
\pi(x)=\int_2^x \frac{dt}{\log t}+\pieter{O\bigl(\sqrt{x}\log x\bigr)},
\end{equation}
where we have given an error term that can be established
if the  Riemann Hypothesis holds true.
By partial integration, \eqref{hadamard} yields an asymptotic
series expansion in the sense of Poincar\'e, with initial term $x/\log x$.
The integral turns out to give a much better numerical approximation to
$\pi(x)$, than does $x/\log x$.
Thus,
perhaps in arriving at Claim \ref{tauclaim} Ramanujan was deceived by
a false analogy with the problem of the distribution of the primes. Some evidence for this is provided by his statement
that the proof of $\sum_{k=1}^n t_k=o(n)$ in the case $q=5$ is ``quite elementary and very similar to that for showing that $\pi(x)=o(x)$'' \cite[p.\,97]{lnbIII}.

 Rankin \cite[p.\,263]{rankintau} was of the opinion that Hardy must have informed Ramanujan of Landau's method after receiving his first
letter and that Ramanujan had a sufficient understanding of the method to be able to make
the informed Claim \ref{tauclaim}.   About his first impressions, Hardy remarked, ``he \dots\ had indeed but the vaguest idea of what a function of a complex variable was \cite[p.~xxx]{cp}.''  However, Hardy further remarked, ``In a few years' time he had a very tolerable knowledge of the theory of functions  \dots \cite[p.~xxxi]{cp}.''  Moreover, on the next-to-last page of Ramanujan's third notebook, which was likely written either shortly before he returned to India or shortly after his arrival home, several integral evaluations are recorded.  Next to one of them appear the words, ``contour integration''  \cite[pp.~391]{nb}.

\subsection{The leading term in Claim \ref{tauclaim}}
\label{sec:revisit}
In this subsection, in greater detail, we consider Claim \ref{tauclaim} for $q\in \{3,5,7,23,691\}$ in
the unpublished manuscript \cite{BerndtOno}.  Still further details can be found in the article by
A.~Ciolan et al.\,\cite[Sec.\,5]{CLM}.
These claims are recorded in
Table \ref{tab:table2}, and all involve the 
tau-function. Not listed are those cases where Ramanujan
only claimed an estimate of the form $O(n\log^{-\delta}n)$.

\begin{table}[ht]
\footnotesize
\renewcommand{\arraystretch}{1.25}
\begin{tabular}{|c|c|c|c|l|c|}\hline
$q$ & $\delta_q$ & E.P.  & $C_q$ & pp. & Sec.\\ \hline \hline
$3$ & $+$ & $+$  &  $+$ 		& 22--23 & 11\\ \hline
$5$ & $+$ & $+$  &   		&  06--08 & 2\\ \hline
$7$ & $+$ & $+$  &  $+$ 		& 11--12 & 6\\ \hline
$23$ & $+$ & $-$  &  $-$ 		&  36--37 & 17\\ \hline
$691$ & $+$ &   &   		&  24--25 & 12\\ \hline
\end{tabular}
\bigskip
\caption{Non-divisibility claims of Ramanujan}
\label{tab:table2}
\end{table}

The $``+''$ entry indicates a correct claim, the $``-''$ entry indicates a false
one, and no entry indicates that no claim was made.
The first column concerns the value of $\delta_q$
(recorded in Table \ref{tab:table3}),
the second the Euler product (E.P.) of the generating series, the third the value of the constant $C_q,$ and the two remaining
columns give the page
numbers, respectively, section numbers in \cite{lnb},  where
the specific claims can be found.

Using his ideas from 
\cite{Rankin61}, Rankin confirmed the correctness
of $C_3$ and $C_7$, and the signs in the $\delta_q$ column
\cite[p.\,10]{Rankin76}. However,
$C_{23}$ needs a minor correction, as first
pointed out by the second author \cite{MRama},
namely, Ramanujan omitted a factor $(1-23^{-s})^{-1}$ in the
generating function (17.6).
He calculated correctly the asymptotic constant associated to his Euler product (17.6), but it has to
be multiplied by $23/22$ in order to obtain
the correct value of $C_{23}$.

  The associated Dirichlet series
$$T_q(s)=\sum_{q\nmid \tau(n)}\frac{1}{n^s}, \quad \Re(s)>1,$$
with
$q\in \{3,7,23\}$, are the easiest in the sense that
the most complicated function they involve is $L(s,\chi_{-q})$, where $\chi_{-q}$  denotes the nontrivial, quadratic Dirichlet character
modulo $q$.
In these three cases, we have $\delta_q=1/2.$
 As $\delta_q=q/(q^2-1)$ for
 $q\not\in \{2,3,5,7,23,691\}$
(see, e.g., Serre \cite[p.\,229]{Ser}),
it follows that for the tau-function there are no
further primes with $\delta_q=1/2$. For $\delta_q<1/2$, the constants
$C_q$ become far more difficult to evaluate, and so this
may explain why Ramanujan did not venture to write down
$C_5$ and $C_{691}.$

\subsection{The error term in Claim \ref{tauclaim}}
In order to deal with Claim \ref{tauclaim} in its sharpest form \eqref{valseanalogie}, we proceed as in our
determination of Shanks' constant in Section \ref{sec:EKsumoftwosquares}.
The \pieter{Dirichlet series} $T_q(s)^{1/\delta_q}$ has a simple pole at
$s=1.$
Employing the product representation for  $T_q(s)^{1/\delta_q}$ and taking logarithmic derivatives, we arrive at a closed expression for the associated Euler--Kronecker constant
(similar to \eqref{theanswer}). It involves terms of the form $\frac{L'}{L}(1,\chi)$, which can be evaluated to a high numerical precision
with the methods in the papers of
Languasco \cite{Languasco2021a}  and A.~Languasco and L.~Righi \cite{LanguascoR2021}.
This then leads to the following result.

\begin{theorem}[Moree \cite{MRama}] For $q\in \{3,5,7,23,691\},$
the set $S_{\tau;q}=\{n:~q\nmid \tau(n)\}$
has an Euler--Kronecker constant
given in Table \ref{tab:table3}.
\end{theorem}
\begin{corollary}
Ramanujan's claim \eqref{valseanalogie} is false if $r>1+\delta_q$.
\end{corollary}
The key to the work described above is that by the congruences
described in Section \ref{taucongruences}, the relevant sets $S_{\tau;q}$ are multiplicative abelian, and hence we
can apply Theorem \ref{mainabelian}.
\begin{table}[ht]
\footnotesize
\renewcommand{\arraystretch}{1.25}
\begin{tabular}{|c|c|c|c|}\hline
set & $\gamma_{S_{\tau;q}}$  & winner &  $\delta_q$\\ \hline \hline
$3\nmid \tau(n)$ & $0.5349(21)\dotsc$ & Landau &  $1/2$\\ \hline
$5\nmid \tau(n)$ & $0.3995(47)\dotsc$ & Ramanujan & $1/4$\\ \hline
$7\nmid \tau(n)$ & $0.2316(40)\dotsc$ & Ramanujan & $1/2$\\ \hline
$23\nmid \tau(n)$ & $0.2166(91)\dotsc$ & Ramanujan & $1/2$\\ \hline
$691\nmid \tau(n)$ & $0.5717(14)\dotsc$ & Landau & $1/690$ \\ \hline
\end{tabular}
\bigskip
\caption{Euler--Kronecker constants related to Claim \ref{tauclaim}}
\label{tab:table3}
\end{table}

The Euler--Kronecker constants in Table \ref{tab:table3} have been confirmed
and determined with slightly higher
accuracy by Ciolan et al.\,\cite{CLM};
further computed digits are indicated in
parentheses.

\section{Fourth Interlude: Squares galore!}
\label{squaresinterlude}
Let $r_k(n)$ denote the number of  representations of $n$ as a sum of $k$ squares, where $k$ is a positive 
 integer. 
Squares of positive numbers, negative
numbers and zero are all allowed, and the ordering of the squares of the numbers
that occur in this summation also counts. 
We have, for example, $r_{24}(2)=1104$.
In Section \ref{circleproblem}, we considered $r_2(n)$.   Focusing, for simplicity, on $r_k(p)$ with $p$ prime, one finds that 
 this is a polynomial in $p$ for several smaller values of 
 $k$. For example, Jacobi in the early part of the 19th century found that $r_4(p)=8p + 8$
 and $r_8(p)=16p^3 + 16$.
However, for larger $k$ there is very often no polynomial formula for $r_k(p)$. Here, Fourier coefficients of cusp forms come to the rescue. A particular charming example, as it involves both sums of squares and Ramanujan's tau-function, is
$$r_{24}(p)=\frac{16}{691}(p^{11}+1)+\frac{33152}{691}\tau(p).$$
By Deligne's bound 
\eqref{Delignebound} for $\tau(p)$, we have $r_{24}(p)=\frac{16}{691}p^{11}+O(p^{\frac{11}{2}})$.
For a very readable elementary introduction on the behavior of $r_{24}(p)$, see B.~Mazur \cite{error}.

The book by E.~Grosswald \cite{Grosswald} and the monograph by S.~Milne \cite{Milne} are perhaps the two primary sources on $r_k(n)$. A shorter read is Chapter 9 of Hardy \cite{Hardy}, which ends with a discussion of 
$r_{24}(n)$.
\section{Some generalizations of Ramanujan's claims}
\label{generale}

Generalizing Ramanujan's Claims \ref{ramaclaim} and \ref{tauclaim}
in Section \ref{sec:LvR}, we presented
the theory of multiplicative sets.
Here we  consider analogues of Claim \ref{ramaclaim} involving the representations of
integers by binary quadratic forms, other than $X^2+Y^2$, and analogues of Claim \ref{tauclaim} for other moduli and/or other Fourier coefficients of modular forms.
Of the many relevant papers, we will discuss
only a modest selection. We start by discussing the relevance of a variant of the original
problem.

\subsection{Easy numerical approximation of \texorpdfstring{$B(x)$}{Bx}} 
Shanks \cite{Shanks} asserted that
``An unsolved problem of interest is to find an approximation to
$B(x)$ that could be computed without undue difficulty by a convergent
process, and which would be accurate to $O(x\log^{-m}x)$ for all $m$.'' Perhaps he would have regarded the following result as giving an adequate answer (where  the function $H(s)$
in the integral
denotes the right hand side in
\eqref{LSsquare}, a function which can be analytically extended up to $\Re s>1/2$, 
so that the integral is well-defined).
\begin{theorem}
\label{zeveneen}
Let $0<\epsilon<1/2$. There exists a constant $c>0$ such that
$$B(x)=\frac{1}{\pi}\int_{1/2+\epsilon}^1 \sqrt{|H(s)|}\,x^{\sigma}\frac{\,d\sigma}{\sigma}+
O(x\,e^{-c\sqrt{\log x}}).$$
If we assume RH for both $\zeta(s)$ and $L(s,\chi_{-4})$, then we can replace the error term by $O(x^{\frac12+\epsilon})$. 
\end{theorem} 
This result can be regarded as a corrected version of
Claim \ref{ramaclaim}! It can be used in practice to 
numerically approximate $B(x)$.
We refer to O.~Gorodetsky 
and B.~Rodgers 
\cite[Appendix B]{GoroRodgers} or  C.~David 
et al.\,\cite{DDNS} for a proof. For a more general
form of it, see G.~Tenenbaum \cite[p.\,291]{Tenenbaum}.
The authors note that Theorem \ref{zeveneen} (known to experts before) 
seems to go back to a 1976 paper of K.~Ramachandra \cite{Ramachandra}.

\subsection{A variant of Claim \ref{ramaclaim}}
Let $d\ge 1$ be a squarefree integer.
By the Hasse principle the \emph{negative Pell equation} $X^2-dY^2=-1$ has a solution with $X$ and $Y$ 
rational if and only if $d$ has no prime factor $p\equiv 3\pmod{4}$. 
A minor variation of Landau's proof of \eqref{Edmund} gives that
the number of such $d\le x$ is asymptotically $\sim cx\log^{-1/2}x$ for some $c>0$. 
P.~Stevenhagen \cite{peterpell} conjectured that a similar asymptotic holds if we ask for \emph{integer} solutions $X$ and $Y$ and indicated an explicit value of $c$. 
Recently P.~Koymans and C.~Pagano \cite{KP} established this
conjecture, after earlier deep work by \'E.~Fouvry and J.~Kl\"uners \cite{FK} placing $c$ in some interval. 

\subsection{Generalizations of Claim \ref{ramaclaim}}
Let $s_1,s_2,\dotsc$ be a sequence of integers that can be written as  sums of two squares. Generalizations mostly
concern other related sequences of integers in the following Sections
\ref{Claim1:quadri} and \ref{Claim1:compu}, or they focus on $s_1,s_2,\dotsc$, but ask for
more refined distributional properties (Section \ref{claim1:distri}). Indeed, the sequence $s_1,s_2,\dotsc$ is among those most intensively studied.

Several of the generalizations discussed here have also been considered in the
\emph{function field} setting; see, for example, the paper of Gorodetsky \cite{poly}.

\subsubsection{Generalizations to other (binary) quadratic forms}
\label{Claim1:quadri}
C.~Stewart and Y.~Xiao \cite{SX} showed that if $F$ is a binary form of degree $d\ge 3$ with integer coefficients and non-zero discriminant, then
$$\{n:|n|\le x,~F(X,Y)=n\text{~for~some~}X,Y\in\mathbb Z\}=C_Fx^{2/d}+O_{F,\epsilon}(x^{\beta_F+\epsilon}),$$
where $C_F$ and the rational number $\beta_F<2/d$ can be obtained explicitly.
It remains to discuss what happens for
arbitrary positive definite binary quadratic forms. Let us first discuss the easier problem
in which integers can be
represented by {\it some} primitive binary
form of a prescribed discriminant $D$. A less difficult variant of this arises
when one restricts oneself
to counting integers that are co-prime only to $D$. In this direction, R.D.~James \cite{James} showed, by Landau's method, that the
number $B_D(x)$ of positive integers $m\le x$, co-prime to $D$, which are
represented by some primitive binary form of discriminant $D\le -3$, satisfies
\begin{equation}
\label{jjames}
B_D(x)=b_D\frac{x}{\sqrt{\log x}}\Bigl(1+O\Bigl(\frac{1}{\sqrt{\log x}}\Bigr)
\Bigr),\quad x\to\infty,
\end{equation}
where $b_D$ is some explicit positive value. The algebraic fact on which the
proof relies is that an integer $m$ that is co-prime to $D$ is represented by a
primitive binary quadratic form of discriminant $D$, if and only if the primes
$p$ that occur with odd exponent in $m$ satisfy $(D/p)=1$, where
$(D/p)$
denotes the Kronecker symbol.  In particular,
$(D/p) $ equals
 the Legendre symbol if  $p$ is an odd prime.
G.\,Pall \cite{Pall2}
showed that (\ref{jjames}) remains valid, but with a different explicit
constant $b_D$, if the restriction that $m$ is co-prime to $D$ is dropped. The
second coefficient in the Poincar\'e series of Pall's theorem was obtained by W.\,Heupel
\cite{Heupel}, who
also obtained the second coefficient in the Poincar\'e series of $B_D(x)$.
K.S.\,Williams
\cite{Williams2} gave
\pieter{a short} proof of \eqref{jjames}
with an error term
$O\bigl((\log \log x)^{-1}\bigr)$.

 Now we return to the problem of  representations by a given binary
quadratic form. In his 1912 Ph.D.\,thesis, P.~Bernays \cite{Bernays} proved that
the counting functions for the integers represented by a reduced binary
form of discriminant $D$ are asymptotically identical; specifically, he showed that they behave
asymptotically as
$$C(D)\frac{x}{\sqrt{\log x}},$$
where $C(D)$ is a positive constant depending only on $D$.
Note that $C(-4)$ is the Landau--Ramanujan constant.
Odoni \cite{O4} presented his own 
proof of Bernays' result in the 
more precise form
\begin{equation}
\label{preciseform}    
C(D)\frac{x}{\sqrt{\log x}}\Big(1+O\Big(\frac{1}{\log^{e(f)}x}\Big)\Big),\quad e(f)>0,
\end{equation}
but it seems that his method does not permit the calculation of the constant $C(D)$ and gives no information about $e(f)$. Fomenko \cite{Fomenko} gave a nice proof of \eqref{preciseform}, which enables one to calculate $C(D)$ and to estimate $e(f)$.

Bernays left number theory for logic (in which he was to become famous) and did not publish his thesis.
Thus,
unfortunately many later researchers were unaware of it, or
did not have access to it. Had they been,
this would have profoundly altered what subsequently happened, as Bernays'
methods are more powerful than those employed by many later authors.
As it is, Bernays' result
was only improved and generalized after roughly 60 years by Odoni, who
wrote a long series of papers generalizing this type of result,
 eventually resulting in his theory
of \emph{Frobenian multiplicative functions}; see 
\cite{Odonishort} and \cite{O1} for a short, respectively, longer survey.

\subsubsection{Generalizations with a computational aspect}
\label{Claim1:compu}
P.~Shiu \cite{Shiu} adapted the Meissel--Lehmer method, initially developed
for calculating $\pi(x)$, to $B(x)$. This allows one to compute $B(x)$
more efficiently, but at the expense of explicitly knowing which integers $\le x$
can be written as a sum of two squares.

In 1966, focusing on the family of forms $X^2+kY^2$, Shanks and L.P.~Schmid \cite{SS} considered the problem of
determining $C(-4k)$ in which  negative
values of $k$ were also considered. They expected that
$C(-8)$ would be the largest such
constant. However, D.~Brink et
al.\,\cite{BMO} showed that
$C(-4k)$ is unbounded as $k$ runs through
both the positive and negative integers.

Let $F_3(x)$ denote the number of integers $n\le x$ represented by $X^2+XY+Y^2$.
Since $C(-4)>C(-3)$, we have $B(x)\ge F_3(x)$ for all
$x$ sufficiently large. In connection with his work on lattices, P.~ Schmutz Schaller \cite{Schmutz} conjectured
that actually $B(x)\ge F_3(x)$ for \emph{every}
$x\ge 1$. This conjecture was proven by
the second author and H.J.J.~te Riele
\cite{MteR}. 
\pieter{An adaptation of the very elementary method of Selberg \cite[pp.\,183-185]{Selberg} alluded to in the introduction allowed them to show only that $B(x)\ge F_3(x)$ for every
$x\ge 10^{9111}$ \cite[\S 10]{MteR}, and so again only got honorable mention.}

\'E.~Fouvry et al.~\cite{FLW} determined an asymptotic for the integers $n\le x$ represented by
\emph{both} $X^2+Y^2$ and $X^2+XY+Y^2$. Languasco and the first author
\cite[Sec.\,7]{LaMoree}
noticed that this, in essence, was already done by Serre \cite[pp.~185--187]{SerreChebotarev} who counted asymptotically
the number of non-zero Fourier coefficients $a_n$ of 
$q\prod_{n\ge 1}(1-q^{12n})^2$. Namely, 
we have $a_n\ne 0$ if and only if $n\equiv 1({\rm mod~}12)$ and $n$ is representable by both $X^2+Y^2$ and $X^2+XY+Y^2$.

Note that $B(x)$ and $F_3(x)$ are related to the quadratic, respectively, hexagonal lattices.
A quantity associated with an $n$-dimensional lattice $L$ is its 
\emph{Erd\H{o}s number}\footnote{Not to be mixed up with the celebrated Erd\H{o}s collaborative distance!} and is given by $E_L=F_L\,d^{1/n}$, where $d$ is the determinant of the lattice and $F_L$ its \emph{population fraction}, which is given by
$$F_L=\lim_{x\rightarrow \infty}\frac{N_L(x)\sqrt{\log x}}{x}\quad\text{if~}n=2,\quad F_L=\lim_{x\rightarrow \infty}\frac{N_L(x)}{x}\quad\text{if~}n\ge 3,$$
where
$N_L(x)$
is the
\emph{population function}
associated with 
the
corresponding
quadratic form, i.e., the number of values not exceeding $x$ taken by the form.
The \emph{Erd\H{o}s number} is the population fraction when the lattice is normalized to
have covolume 1.
In the case of a two dimensional lattice, it
is related to the constant $b_D$ appearing in \eqref{jjames} (see e.g.,
\cite{conway}). In 2006, the second author and
R.~Osburn \cite{MO} showed, relying on the
work of many others, that
the Erd\H{o}s number is minimal for the hexagonal lattice. Also, relying on the work of many others, J.H.~Conway and N.J.A.~Sloane \cite{conway} solved the analogous  problem for dimensions 3 to 8 in 1991.

\subsubsection{Other distributional aspects of sums of two squares}
\label{claim1:distri}
A pair $(n,n+1)$ is said to be \emph{$B$-twin} if both $n$ and $n+1$ can be written as a sum  of two squares. In 1965, G.J.~Rieger \cite{Rieger} proved with the one-dimensional sieve the upper bound
$\ll x/\log x$ for the number of such pairs with $n\le x$. The analogous lower bound was independently proved in 1974 by C.~Hooley \cite{Hooley3} using the asymptotic estimate 
$$\sum_{n\le x}r(n)r(n+1)=8x+O(x^{5/6+\epsilon})$$
due to T.~Estermann \cite{Estermann}, and by K.-H.~Indlekofer \cite{Indlekofer} with the sieve method. It is a very challenging problem to prove an asymptotic estimate for
$r(n)r(n+h)r(n+k)$ with $n\le x$ and 
$h\ne k$ fixed positive integers. Hooley \cite{Hooley2} made some progress by showing that this 
function, for fixed $h$ and $k$, is infinitely often positive (thus answering a 
question of J.E.~Littlewood).


Given any two positive integers $k$ and $l$, we let
$B(x,k,l)$ denote the number
of integers $s_i\le x$ with
$s_i\equiv l\,({\rm mod~}k)$.
It is sufficient to assume that $k$ and $l$ are co-prime. We will assume that $l\equiv 1\,({\rm mod~}4)$ in the case $k\equiv 0\,({\rm mod~}4)$,
for otherwise $B(x,k,l)=0$ for every $x$.
K.~Prachar \cite{Prachar} showed that
\begin{equation}
    \label{uniform}
B(x,k,l)\sim B_k\frac{x}{\sqrt{\log x}},\quad x\to\infty,
\end{equation}
where $B_k$ is a positive constant depending
only on $k$. His Ph.D.~student, H.~Beki\'c,
\cite{bek} proved that this asymptotic holds uniformly in the range $k$ up to $e^{c\sqrt{\log x}}$, where $c$ is a positive constant, a result subsequently extended by Prachar \cite{Prachar2} to $e^{(\log x)^{2/3-\epsilon}}$. Iwaniec \cite{Iwaniec}, using the half dimensional sieve, strengthened this to
$$B(x,k,l)=B_k\frac{x}{\sqrt{\log x}}\Bigl(1+O\Bigl(\Bigl(\frac{\log k}{\log x}\Bigr)^{1/5}\Bigr), \quad x\to\infty,$$
where the implicit constant is absolute.

Gorodetsky \cite{Goro2} showed that $B(n,k,l_1)>B(n,k,l_2)$, with $l_1,l_2$ a quadratic, respectively,
non-quadratic residue modulo $k$, for a density-1 set of integers $n$, provided some standard number theoretic conjectures including GRH hold true.

J.~Schlitt \cite{Schlitt} generalized \eqref{uniform} to other positive definite binary quadratic forms.
However, although asymptotically there is equidistribution, at the level of the second order
constant, this is not always the case.
For example, if $k\equiv 1\pmod{4}$ is a prime, then there is a preponderance
for residue classes $0\pmod{k}$ over the others. Numerically, this effect is clearly detectable.

Let $r_p(n)$ be the number of representations of $n$ as a sum of two prime
squares. It is known that asymptotically
$\sum_{n\le x}r_p(n)^j\sim 2^{j-1}\pi x\log^{-2}x$ for $j=1,2,3$,
where the case $j\le 2$ is fairly standard and $j=3$ is due to
V.~Blomer and J.~Br\"udern \cite{BB}. Recently, partial results for $j\ge 4$ have
been obtained by C.~Sabuncu \cite{Sabuncu}.
A.~Sedunova \cite{Sedunova22} showed that asymptotically
$\sum_{n\le x}r(n)r_p(n)\sim  12Gx/(\pi^2\log x)$, with $G$ the Catalan constant. Together with Granville and Sabuncu \cite{GSS}, she gave a precise estimate for the number of integers $\le x$ that can written as the sum of a square and the square of a prime. This estimate involves the multiplication table constant that arises in studies of the number of distinct integers in the $N$-by-$N$ multiplication table.

P.~L\'evy \cite{Levy} gave a simple heuristic derivation of 
\eqref{Edmund}, however without determing $K$. His argument also led him to conjecture that $R_k(x)$, the number of $n\le x$ for which $r_2(n)=k$, asymptotically satisfies
$$R_k(x)\sim \frac{Kx}{\sqrt{\log x}}\frac{e^{-\theta}\,\theta^k}{k!},\quad\text{~where~}\theta=c\sqrt{\log x}.$$
In probabilistic terms, this means roughly that the integers $m$ for which $r_2(m)$ has a specified value, have a Poisson distribution with parameter $\theta$. W.J.~LeVeque \cite{LeVeque} showed that the conjecture is wrong and that in fact the asymptotic behavior of $R_k(x)$ not only depends on the size of $k$,
but also its arithmetic structure.

W.D.~Banks et al.\,\cite{BHMN} showed that every integer $n>720$ that can be written
as a sum of two squares satisfies Robin's inequality
\begin{equation}
\label{Robin}
\sum_{d|n}\frac{1}{d}<e^{\gamma}\log \log n.
\end{equation}
G.~Robin \cite{Robin} famously proved that
the Riemann Hypothesis is true if and only if  \eqref{Robin} holds for all
$n>5040$. We leave it as a (not so difficult) challenge for the reader to show
that \eqref{Robin} is satisfied for all \emph{odd} integers, and thus, to wit, prove
half of the Riemann Hypothesis!

A.~Balog and T.D.~Wooley \cite{BW}
found “unexpected irregularities” in the distribution of
the \pieter{integers} $s_i$ in short intervals. 
To be precise, they showed that there are infinitely many short
intervals containing considerably more
integers $s_i$ than expected, and infinitely many intervals containing
considerably fewer than expected. \pieter{They used
a method of H.~Maier \cite{Maier}, which imposes restrictions on the equidistribution of primes, and has been since then frequently used (``Maier's matrix method''). For example, Granville and 
Soundararajan \cite{GS} used it to establish an ``uncertainty principle'' saying that most arithmetic sequences of interest are either not-so-well distributed in longish arithmetic progressions, or are not-so-well distributed in both short intervals and short arithmetic progressions.}
Hooley \cite{Hooley1,Hooley2,Hooley3,Hooley4} wrote four papers on  the distribution of the gaps
$s_{i+1}-s_i$.

C.~David et al.~\cite{DDNS} studied how frequently
$(s_i,s_{i+1})$ belongs to a prescribed congruence class pair modulo $q$.
Certain congruence pairs appear more frequently than others, which they heuristically 
explain. A similar
phenomenon was found earlier by R.J.~Lemke Oliver and K.~Soundararajan for consecutive prime numbers \cite{LOSound} and caused quite a sensation.
N.~Kimmel and V.~Kuperberg \cite{KK} considered the situation where given a fixed integer $m\ge 0$, one requires $s_i,s_{i+1},\ldots,s_{i+m}$ to be in prescribed progressions (which are allowed to be different).

The energy levels of generic, integrable systems are conjectured to be Poisson distributed in
the semiclassical limit.
The square billiard, although completely integrable, is non-generic in this respect.
Its levels, when suitably scaled, are the numbers $s_i$. This led to some interest in theoretical mathematical physics in the distribution of the numbers $s_i$ \cite{CK1,CK2,FKR}.
For example, T.~Freiberg
et al.\,\cite{FKR} formulate an analog of the Hardy--Littlewood prime $k$-tuple conjecture for sums of two squares, and show that it implies that the spectral gaps, after removing degeneracies and rescaling, are Poisson distributed.

Let $s(n)$ denote the sum of the proper divisors of $n$ (thus $s(n)=\sigma_1(n)-n$). 
Its study goes back a long time,
beginning with the investigations of the ancients into perfect numbers and amicable pairs. P. Erdős et al.\,\cite{EGS} conjectured that if 
$A$ is a subset of the natural numbers of asymptotic density zero, then the set of integers $m$ for which $s(m)\in A$ also has asymptotic density zero. Troupe \cite{Troupe} proved this conjecture in case $A$ is 
the set of integers which can be represented as the sum of two squares.

A positive integer $n$ is called \emph{square-full} if $p^2$
divides $n$, whenever $p$ is a prime divisor of $n$. It was shown in 1935
by P.~Erd\H{o}s and G.~Szekeres \cite{ES} that the number of square-full integers $\le x$
is equal to $\zeta(3/2)\sqrt{x}/\zeta(3 )+O(x^{1/3})$. Thus the square-full
numbers are not much more abundant than perfect squares. Thus one might
conjecture,
as did Erd\H{o}s, that $V(x)$ the number of integers
 $\le x$ that can be written as a sum of two square-full integers grows
asymptotically like $cx/\sqrt{\log x}$. However, V.~Blomer \cite{BII} showed that 
asymptotically 
$$V(x)=x\log^{-\alpha+o(1)}x\text{~~with~~}\alpha=1-2^{-1/3}=0.2062994740\ldots.$$
It is a famous result of C.-F.~Gauss 
that $n$ can be written as a sum of three integer squares if and only if $n\ne 4^a(8b+7)$. 
D.R.~Heath-Brown \cite{H-B} showed that there is an
effectively computable constant 
$n_0$ such that every $N\ge n_0$ is a sum of at most three square-full integers.

In 1922, Hardy and Littlewood
\cite{HL} in a now famous paper, conjectured an asymptotic formula, suggested by a formal application of their circle method, for the number of integers of the
form $p+m^2+n^2$, with $p$ prime. Under GRH, Hooley \cite{Hooley57}, see also his
book \cite[Chp.\,5]{Hooleybook}, confirmed the validity of this formula
and several years later Linnik \cite{Linnik} did so unconditionally. Very recently Heath-Brown \cite{H-B2} examined the distribution of the solutions $p,m,n$.

\subsection{Generalizations of Claim \ref{tauclaim}}
\label{sec:claim2general}
In this section we assume some familiarity with the theory
of modular forms, the reader can consult, for example, the 
books by F.~Diamond and J.~Shurman \cite{DS}, N.~Koblitz \cite{Koblitz} or the relatively short classic book by Serre \cite{course}. 
More advanced is the book by Ono \cite{Onoweb}, which has some focus on Fourier coefficient congruences.
Ram Murty \cite{Murtysieving}  gives a nice warm-up for what is going on in this section, where both the algebraic and analytic side are kept at a very accessible level.

\subsubsection{Non-divisibility of sums of divisors functions}
\label{divi}
For the primes $q=2,3,5,7, 691$, the Ramanujan congruences
\eqref{ramacongruences} relate
the non-divisibility of certain Fourier coefficients
to those of
$n^a\sigma_k(n)$ for appropriate $a$ and $k$. In essence, we have an Eisenstein series of weight $k+1$ (having $\sigma_k(n)$ as its $n$th Fourier coefficient) that modulo an appropriate prime becomes a cusp form.
Thus we are lead to consider, for an \emph{arbitrary} integer $k\ge 1$ and a prime $q,$
\begin{equation}
\label{Skn}
S_{k,q}(x) \coloneq S_{\sigma_k;q}(x)=\sum_{\substack{n\le x,\\q\nmid \sigma_k(n)}}1.
\end{equation}
By a minor variation,  the general case
$q\nmid n^a\sigma_k(n)$ can be also handled.
It is therefore natural that Ramanujan was interested in the
asymptotic behavior of
$S_{k,q}(x)$, and he seems to have been the first to do so.
He discusses this function
in his unpublished manuscript
\cite[Sec.\,19]{BerndtOno},
and made three claims (also reproduced by
Rankin \cite{Rankin76}), which were proved
in 1935 by
Watson \cite{watson1}. One of these claims asserted that, for odd $k$,
\begin{equation}
\label{skestimate}
S_{k,q}(x)=O(x\log^{-1/(q-1)}x),
\end{equation}
and it is discussed by Hardy in his Harvard lectures on Ramanujan's mathematics \cite[\S 10.6]{Hardy}.

 The precise asymptotic behavior of
$S_{k,q}(x)$ was first determined by Rankin \cite{Rankin61}.
His Ph.D.~student, Eira Scourfield, \cite{Scourfield64} generalized his work by establishing an asymptotic formula in the case where
a prime power
exactly divides
$\sigma_k(n)$.

Put $h=\frac{q-1}{(q-1,k)}$. Ciolan et al.\,\cite{CLM} went beyond determining an asymptotic formula, and showed that when $h$ is even and  $x\to\infty$,
\begin{equation}
\label{starrie2}
S_{k,q}(x)=\frac{c_0\,x}{\log^{1/h}x}\Bigl(1+\frac{1-\gamma_{k,q}}{h\log x}+\frac{c_2}{\log^2 x}+\dotsm+
\frac{c_j}{\log^j x}+O_{j,k,q}\Bigl(\frac{1}{\log^{j+1}x}\Bigr)\Bigr),
\end{equation}
where $\gamma_{k,q}$ is the Euler--Kronecker constant for $\sigma_k(n)$, and can be explicitly given. The case when $h$ is odd is rather trivial, and there we have
$S_{k,q}(x)\sim c_{k,q}x$, where $c_{k,q}$ is a positive constant \cite[Sect.\,3.8]{CLM}.

\subsubsection{Analogues of 
Claim \ref{tauclaim} for cusp forms of higher weight for the full modular 
group}
The space of cusp forms of weight
$k$ of the full modular group  is one-dimensional if and only if  $k\in\{12,16,18,20,22,26\}$. For these weights the
`elementary' congruences of prime modulus were 
completely classified using $\ell$-adic representations by the efforts of Serre and Swinnerton-Dyer, cf.~\cite{S-D}.
Ciolan et al.\,\cite{CLM}, using the results on the non-divisibility of the the sum of divisor functions alluded to in $\S \ref{divi}$, determined the associated leading constant and Euler--Kronecker constant (cf.\,\eqref{initstarrie}) for all of these, with one   
exception. The exception is the congruence (rather, congruential restriction) 
\begin{equation}
\label{Haberlandcong}    
a_p(\Delta E_4)^2\equiv 0,p^{15},2p^{15},4p^{15} \pmod{59} \quad (p\ne 59), 
\end{equation}
where $\Delta E_4$ is the unique normalized cusp of weight 16 for the full modular group, which was recently dealt with by S.~Charlton et al.\,\cite{CMM}.
In this case the associated generating series turns
out to be expressible in terms of Dedekind zeta functions of some non-abelian number fields, and one cannot just do  with Dirichlet $L$-series, as in the other cases\footnote{Recall that the Dedekind zeta function of an abelian number field factorizes in Dirichlet $L$-series; for a non-abelian number field Artin $L$-series arise as factors.}. In order for the associated Euler--Kronecker constant of 
$S=\{n:59\nmid a_n(\Delta E_4)\}$ to be obtained with moderate precision, the authors had to assume the Riemann Hypothesis for the fields involved and used a method of Y.~Ihara \cite{Ihara}. 

\subsubsection{Parity of Fourier coefficients and the partition function}
\label{sec:parity}
Let us consider the parity of $\tau(n)$ first.
Trivially $(1-q^n)^8\equiv 1-q^{8n} \pmod{2}$.
By the Jacobi Triple Product Identity (see \cite[Thm.\,2.8]{Andrews}), 
$$\eta(8z)^3=q\prod_{j=1}^{\infty}(1-q^{8j})^{3}=\sum_{k\ge 0}(-1)^k(2k+1)q^{(2k+1)^2},$$
which one can then use to deduce that
\begin{equation}
\label{Deltaparity}    
\Delta(z)\equiv \eta(8z)^3\equiv \sum_{k\ge 0}q^{(2k+1)^2}\pmod{2}.
\end{equation}
It follows that $\tau(n)$ is odd if and only if $n$ is an odd square. We will now
see that this behavior is quite exceptional among
modular forms $f$ of level one with integer Fourier coefficients $a_n$.
J.~Bella\"iche and J.-L.~Nicolas  \cite{BN}, improving on 
an earlier result of Serre \cite[\S 6.6]{Ser},
prove that for such $f$ one
has the asymptotic estimate $$\#\{n\le x:2\nmid a_n\}=C_f\frac{x}{\log x}(\log \log x)^{g(f)-2}\Big(1+O\Big(\frac{1}{\log \log x}\Big)\Big),$$
if $g(f)\ge 2$. 
Here $g(f)$ is the order of nilpotency of $f$ (the 
minimum number of Hecke operators that are guaranteed to make $f$ vanish, no matter which
ones you pick).
They show that $g(f)=1$ only if $f=\Delta$, which by \eqref{Deltaparity} leads to
$\#\{n\le x:2\nmid \tau(n)\}=\sqrt{x}/2+O(1)$.
They first establish their asymptotic estimate for $f=\Delta^k$ with $k$ arbitrary and odd.
They then try to write a general $f$ (taken modulo two) in a favorable way as linear combinations
of powers of $\Delta$, and in this way obtain the general result. They use here the simple result that the
graded algebra of modular forms modulo $2$ on $\text{SL}_2(\mathbb Z)$ is ${\mathbb F}_2[\Delta]$ (see
Swinnerton-Dyer \cite[Thm.\,3]{S-D-ladic}).
For some further remarks see, for example, Ono \cite[\S 2.7]{Onoweb}.

Although it is a bit tangential, we will briefly discuss the parity of the partition function $p(n)$, since it is attracting a lot of research interest, and 
much remains to be done.
The parity
of $p(n)$ seems to be quite random, and it is widely believed that the partition function is “equally
often” even and odd. More precisely, T.R.~Parkin and Shanks \cite{PS} 
made the conjecture that
$$\#\{n\le x:2\nmid p(n)\}\sim \frac{x}{2}.$$
This is far from being proved: at the moment of writing it cannot be excluded that there exists an
$\epsilon>0$ such that the latter counting function 
is $O(x^{1/2+\epsilon})$ for some $\epsilon>0$.
Bella\"iche and Nicolas \cite{BN}, using their approach involving powers of $\Delta$, showed
that, for $x\ge 2$,
$$\#\{n\le x:2\mid p(n)\}\ge 0.069\sqrt{x}\log \log x,\quad \#\{n\le x:2\nmid p(n)\}\ge \frac{0.048\sqrt{x}}{\log^{7/8}x}.$$
The latter result was considerably improved by Bella\"iche et
al.\,\cite{BGS} who showed that the latter right hand side can be replaced by $\gg \sqrt{x}/\log \log x$ and that this also holds for the number of Fourier coefficients of any weakly holomorphic modular form of half-integral weight coming from the ring of integers of a number field with $2$ replaced by any prime ideal.

S.~Radu \cite{Radu} proved that every arithmetic progression $r\pmod{t}$
contains
inﬁnitely many integers $N$ for which $p(N)$ is even, and inﬁnitely
for which it is odd (this was a conjecture of M.~Subbarao \cite{subbarao}).
It is an open problem to determine an upper bound for the smallest such $N$.

\subsubsection{Non-divisibility of integer valued multiplicative functions}
Analogues of Claim \ref{tauclaim} for $S=\{n:q\nmid f(n)\}$ can be made, where $q$ is any prime and $f$ is any multiplicative function. In Example \ref{Fordexample} we already discussed the case where $f$ is the Euler totient function and $q$ is any prime. Scourfield \cite{Scourfield77} considered the
case where $f(n)=r_k(n)$, with $r_k(n)$ defined as in Section \ref{squaresinterlude} and the more general case in \cite{Scourfield72}.
In case $f(p)$ is a polynomial the relevant generating series factorizes in terms of Dirichlet $L$-series. A more interesting class is that of the Frobenian multiplicative functions, where $f$ is Frobenian (relative to an extension $K$ of the rationals). This entails that for all primes $p$ and $q$ co-prime to some prescribed number that have the same Frobenius symbol we have $f(p^n)=f(q^n)$ for every $n\ge 1$. 
Now the relevant generating series will factorize as a product of Artin $L$-series.
The density of primes with a prescribed Frobenius symbol is given by the Chebotarev density theorem, see 
P.~Stevenhagen and H.W.~Lenstra \cite{SL} for an introduction.

\subsubsection{Lacunarity} 
If 
$f$ is any multiplicative function, then clearly the 
set $S=\{n~:~a_n\ne 0\}$ is multiplicative.
The set $S$ is lacunary if it has natural density zero. 
Serre \cite{serrelacunary} classified all lacunary even powers 
$\eta(z)^r$ of the eta-function with $r>0$. He showed that 
$\eta(mz)^r$ (with $m$ chosen such that only integer powers of $q$ appear in the Fourier series) are lacunary if and only if $r\in \{2,4,6,8,10,14,26\}$.
In the middle of Section \ref{Claim1:compu} we mentioned
$\eta(12z)^2$ (so the case $r=2$). The equivalence given there immediately shows it is lacunary. Serre also writes down asymptotics for the associated sets $S$ of indices of non-zero Fourier coefficients, except for
$r=26$, where determining the asymptotics is an open problem to this day.
In that case Serre showed that $$c_1\frac{x}{\sqrt{\log x}}\le S(x)\le c_2\frac{x}{\sqrt{\log x}},$$ for every $x\ge 2$ and some constants $0<c_1<c_2$.

More generally, if $f=\sum_{n=0}^{\infty}a_nq^n$ is a holomorphic modular form of integral weight $k\ge 0$ for some level $\Gamma_1(N)$ and with integer coefficients $a_n$, then the sequence $a_n\bmod{p}$ is lacunary. This was established by Serre \cite{Ser}, who gave the upper bound 
\begin{equation}
    |\{n\le x: p\nmid a_n\}|\ll \frac{x}{(\log x)^{\beta}},
\end{equation}
where $\beta>0$ is a constant that may depend on $f$. 
Bella\"iche and Soundararajan \cite{Bellsound} considerably
improved this by showing that there exists a rational number $\alpha(f)\in (0,3/4]$, an integer $h(f)\ge 0$, and a positive real constant $c(f)>0$ such that
\begin{equation}
    |\{n\le x: p\nmid a_n\}|\sim c(f)\frac{x}{(\log x)^{\alpha(f)}}(\log \log x)^{h(f)},
\end{equation}
generalizing a result of Serre \cite{Ser} who proved this with $h(f)=0$ in the case when $f$ is an eigenform of all Hecke operators 
$T_m$.

\subsection{Averaging multiplicative functions}
\label{amf}
In Theorem \ref{een} we considered the counting function $S(x)=\sum_{n\le x}i_S(n)$, where $i_S$ assumes non-negative values only and is multiplicative. Likewise, we can ask if a similar result can be obtained if we replace $i_S$ by any non-negative real valued multiplicative function.

Most authors who have been working in this direction,
e.g., B.M.~Bredikhin, H.~Delange,
A.S.~Fainleib,
H.~Halberstam, B.V.~Levin, R.W.K.~Odoni, 
J.M.~Song and
E.~Wirsing, have concentrated on finding conditions on
$\sum_{p\le x}f(p)$ or 
$\sum_{p\le x}f(p)/p$ that are as weak as possible, so  that they could prove an asymptotic formula for
$\sum_{n\le x}f(n).$  In particular, see the book by A.G.~Postnikov \cite{P}. The next result is a famous example of this, and is
due to Wirsing \cite{Wirsing}, \pieter{see also \cite[pp.\,195--215]{P} for a full proof.}

\begin{theorem}
\label{wirsi}
Let $f(n)$ be a multiplicative function such that $f(n)\ge 0$, for
$n\ge 1.$  Suppose that there exist constants $\gamma_1$ and $\gamma_2,$ with
$\gamma_2<2,$
 such that for
every prime $p$ and every $\nu\ge 2,$ $f(p^{\nu})\le
\gamma_1\gamma_2^{\nu}.$
Assume that, as $x\rightarrow \infty,$
\begin{equation}
\label{wirsicondition}    
\sum_{p\le x}f(p)\sim \delta\frac{x}{\log x},
\end{equation}
where $\delta$ is a positive constant.  Then, as $x\to\infty$,
$$\sum_{n\le x}f(n)\sim \frac{e^{-\gamma \,\delta}}{\Gamma(\tau)}\,
\frac{x}{\log x}\,\prod_{p\le x}\Bigl(1+\frac{f(p)}{p}+
\frac{f(p^2)}{p^2}+\frac{f(p^3)}{p^3}+\dotsm  \Bigr).$$
\end{theorem}
Using a theorem of F.~Mertens, see
for example \cite[p.~466, Theorem 429]{hw}, 
$$\prod_{p\le x}\bigl(1-\frac{1}{p}\bigr)\sim \frac{e^{-\gamma}}{\log x},\quad x \to \infty,$$
we see that this result implies Theorem \ref{een}.

In a follow-up paper Wirsing \cite{Wirsing2} proved a variant of Theorem \ref{wirsi}, where the PNT type condition \eqref{wirsicondition} is replaced by the weaker Mertens type condition
$$\sum_{p\le x}\frac{\log p}{p}f(p)\sim \delta\log x.$$
G.~Tenenbaum \cite{Tenenbaumeffective} was the first to establish an effective version of 
this result.
In case $f(n)\in \{0,1\}$ for every $n$, Wirsing's variant is actually stronger than Theorem \ref{wirsi}.
Another variant was established by Song \cite{Song} (Theorem A of part I).
She also extended this type of result to the counting of 
integers having largest prime factor $\le y$ (so-called \emph{$y$-friable integers}), 
\pieter{see also \cite{HTW,TW}.}
Very recently the case where $f$ is real and oscillatory, respectively complex has been also studied \cite{BT,BT2}.

In many situations
one can prove a much more precise estimate for $\sum_{p\le x}f(p)$ (or a related counting function), and this can be used to obtain
a much stronger result than that given in the conclusion of Wirsing's theorem. For example, A.~Granville and D.~Koukoulopoulos \cite{GK} consider the case where $\sum_{p\le x}f(p)\log p=\alpha x+O(x\log^{-A}x)$, for some $A>0$ (their results were improved 
by de la Bret\`eche and Tenenbaum \cite{BT3}). 
There are several methods that can be used in this setting.  One can be regarded as
a considerable extension of the method of contour integration that Landau
used in his paper \cite{L} in 1908, and is commonly called
the Selberg--Delange method.
It applies, in principle, to any arithmetic function whose Dirichlet series is close to a complex power of zeta, 
but requires the analytic continuation of this series,
a condition that can be difficult to confirm in practice.
In view of our story up to this point, the reader might
be inclined to think that Landau--Selberg--Delange (LSD) is a more appropriate name, and
is indeed not alone in this; see the article by A.~Granville and D.~Koukoulopoulos \cite{GK}. 
However, there are also experts arguing
that the name Landau should not be included. 

For introductory material on the
\pieter{(Landau)}--Selberg--Delange method we refer to the books by Koukoulopoulos \cite{Koukou} and Tenenbaum \cite{Tenenbaum}. For a brief exposition, see the article by Finch et al.\,\cite{FMS}.\\

\noindent \noindent {\footnotesize {\it Acknowledgment}.
An important source for this article is 
\cite{MC}, written a quarter of a century ago by the second author and Cazaran. Many things happened in this subject area since then!
We thank Alessandro Languasco for updating us on the high precision computation of the constants of Cilleruelo, Landau--Ramanujan and Shanks and lots of further help. Section
\ref{sec:EKsumoftwosquares} was written in intensive exchange with him.
Further, we thank Maarten Derickx for some helpful correspondence.
The TeX orginal of Table \ref{tab:table1} was kindly sent us by Bernhard Heim. 
Will Craig updated us on recent work on excluded $\tau$-values and 
Lukas Mauth on the parity of the partition function.
Zeev Rudnick gave helpful feedback on \S \ref{sec:Cilleruelo}.
Mike Daas, Nikos Diamantis, Harold Diamond, Ofir Gorodetsky, Anna Medvedovsky, Oana Padurariu, Cihan Sabuncu and Don Zagier commented on an earlier version. Thanks to Ofir, the articles 
\cite{Fomenko,O4,O,Song,Ramachandra,Wirsing2} are now also mentioned. 
In person meetings
of the second author with J\'anos Pintz and Efthymios Sofos  
also led to the addition of interesting further material. The authors are grateful to G\'erald Tenenbaum for helping to bring
Section \ref{amf} in a shape more palatable to 
Selberg--Delange method experts.
The MPIM librarian Anu Hirvonen kindly helped with getting hold of several papers in the references.
The initial impetus for writing this 
article came from Ken Ono.

\begin{thebibliography}{999}
\bibitem{AT} A.~Akbary and Y.~Totani, Binary quadratic forms of odd class number, \url{https://arxiv.org/abs/2408.00184},  preprint.

\bibitem{Alkan}
E.~Alkan, On the sizes of gaps in the Fourier expansion of modular forms, \textit{Canad. J. Math.} \textbf{57} (2005), 449--470.

\bibitem{Alkan2}
E.~Alkan, Average size of gaps in the Fourier expansion of modular forms, \textit{Int. J. Number Theory} \textbf{3} (2007), 207--215.

\bibitem{AZ} E.~Alkan and A.~Zaharescu, 
Nonvanishing of the Ramanujan tau function in short intervals, 
\textit{Int. J. Number Theory} \textbf{1} (2005), 45--51.


\bibitem{AB} 
G.~Almkvist and
B.C.~Berndt, Gauss, Landen, Ramanujan, the arithmetic-geometric mean, ellipses, $\pi$, and the Ladies diary,
\textit{Amer. Math. Monthly}
\textbf{95} (1988), 585--608.

\bibitem{Andrews}
G.E.~Andrews, \textit{The Theory of Partitions}, Encyclopedia of Mathematics and its Applications \textbf{2}, Addison-Wesley Publishing Co., Reading, Mass.-London-Amsterdam, 1976. 


\bibitem{lnbIII}
G.E.~Andrews and B.C.~Berndt, \textit{Ramanujan's Lost Notebook, Part III},
 Springer, New York, 2012.

\bibitem{ayoub}
R.~Ayoub, The lemniscate and Fagnano's conributions to elliptic integrals, \textit{Arch.~Hist.~Exact Sci.} \textbf{29} (1984), 131--149.

\bibitem{BalogOno} A.~Balog and K.~Ono, 
The Chebotarev density theorem in short intervals and some questions of Serre, \textit{J. Number Theory} 
\textbf{91} (2001), 356--371.


\bibitem{BW} A.~Balog and T.D.~Wooley,
Sums of two squares in short intervals, \textit{Canad.~J.~Math.}\,
\textbf{52} (2000), 673--694.



\bibitem{BCO} J.S.~Balakrishnan, W.~Craig and K.~Ono, Variations
of Lehmer's conjecture for Ramanujan's tau-function,
\textit{J. Number Theory} \textbf{237} (2022), 3--14.

\bibitem{BCO2} J.S.~Balakrishnan, W.~Craig, K.~Ono and W.-L. Tsai, Variants of Lehmer’s speculation for newforms, \textit{Adv. Math.} \textbf{428} (2023), 109141.

\bibitem{BOT} J.S.~Balakrishnan, K.~Ono and W.L.~Tsai, Even values of Ramanujan's tau-function,
\textit{Matematica} \textbf{1} (2022), 395--403.

\bibitem{BHMN} W.D.~Banks, D.N.~Hart, P.~Moree and C.W.~Nevans, The
Nicolas and Robin inequalities with sums of two squares,
\textit{Monatsh. Math.}\, \textbf{157} (2009), 303--322.

\bibitem{bek} H.~Beki\'c, \"Uber die Anzahl den Zahlen $a^2+b^2$
in einer arithmetischen Reihe grosser Differenz, 
\textit{J.~Reine Angew.~Math.}\, \textbf{ 226} (1967), 120--131.

\bibitem{BGS} 
J.~Bella\"iche, B.~Green and K.~Soundararajan, non-zero coefficients of half-integral weight modular forms 
mod $\ell$, \textit{Res. Math. Sci.} 
\textbf{5} (2018), no. 1, Paper No. 6, 10 pp.


\bibitem{BN} 
J.~Bella\"iche and J.-L.~Nicolas, Parit\'e des coefficients de formes 
modulaires, \textit{Ramanujan J.} \textbf{40} (2016), 1--44.

\bibitem{Bellsound} 
J.~Bella\"iche
and K.~Soundararajan, The number of non-zero coefficients of modular forms (mod $p$),
\textit{Algebra Number Theory} \textbf{9} (2015), 1825--1856.



\bibitem{oddies} M.A.~Bennett, A.~Gherga, V.~Patel and S.~Siksek, 
Odd values of the Ramanujan tau function, 
\textit{Math. Ann.} \textbf{382} (2022), 203--238.

\bibitem{Bernays} P.~Bernays, \"Uber die Darstellung von positiven,
ganzen Zahlen durch die primitiven, bin\"aren quadratischen Formen einer
nicht-quadratischen Diskriminante, Ph.D.~Thesis, G\"ottingen, 1912, 

\bibitem{Berndt} B.C.~Berndt, \textit{Ramanujan's Notebooks.} IV,
Springer-Verlag, New York, 1994.


\bibitem{bhargava1}
B.C.~Berndt and S.~Bhargava, Ramanujan's inversion formulas for the lemniscate and allied
functions, \textit{J.~Math.~Anal.~Applics.}~{\bf 160} (1991),
504--524.

\bibitem{bkz1}
B.C.~Berndt, S.~Kim, and A.~Zaharescu, Weighted divisor sums and Bessel function series, II, \textit{Adv.~Math}.~\textbf{229} (2012),
2055--2097.

\bibitem{bkz2}
B.C.~Berndt, S.~Kim, and A.~Zaharescu,
Weighted divisor sums and Bessel function series,
IV, \textit{Ramanujan J}.~\textbf{29} (2012), 79--102.

\bibitem{bkz3}
B.C.~Berndt, S.~Kim, and A.~Zaharescu,
 Circle and divisor problems, and double series of Bessel
functions, \textit{}{Adv.~Math}.~\textbf{236} (2013),
24--59.

\bibitem{bkz4}
B.C.~Berndt, S.~Kim, and A.~Zaharescu,
 Weighted divisor sums and Bessel function series,
III, \textit{J.~Reine Angew.~Math}.~\textbf{683}
(2013), 67--96.


\bibitem{bkz5}
B.C.~Berndt, S.~Kim, and A.~Zaharescu,
Weighted divisor sums and Bessel Function Series, V, \textit{J.~Approx.~Thy}.~\textbf{197} (2015), 101--114.

\bibitem{circledivisorsurvey1}
B.C.~Berndt, S.~Kim, and A.~Zaharescu, The circle problem of Gauss and the divisor problem of Dirichlet--still unsolved,
\textit{Amer.~Math.~Monthly} \textbf{125} (2018), 99--114.

\bibitem{BerndtOno} B.C.~Berndt and K.~Ono, Ramanujan's unpublished manuscript on the partition and tau functions with proofs and commentary. The Andrews Festschrift (Maratea, 1998), 
\textit{S\'em. Lothar. Combin.} \textbf{42} (1999), Art. B42c, 63 pp., available at 
\url{https://www.mat.univie.ac.at/~slc/wpapers/s42berndt.pdf}.

\bibitem{Letters} B.C.~Berndt and R.A.~Rankin,
\textit{Ramanujan. Letters and commentary}, History of Mathematics \textbf{9},
AMS, Providence, RI, 1995.



\bibitem{BII} V.~Blomer, Binary quadratic forms with large discriminants and sums of two squareful numbers. II, 
\textit{J.~London Math.~Soc.}\,(2) \textbf{71} (2005), 69--84.

\bibitem{BB} V.\,Blomer and J.\,Br\"udern,
Prime paucity for sums of two squares,
\textit{Bull. Lond. Math. Soc.} \textbf{ 40} (2008), 457--462



\bibitem{BoBo} 
J.M.~Borwein and P.B.~Borwein,
 \textit{Pi and the AGM. A study in analytic number theory and computational complexity}. Canadian Mathematical Society Series of Monographs and Advanced Texts, A Wiley-Interscience Publication,
 John Wiley \& Sons, Inc., New York, 1987.

\bibitem{BMO} D.~Brink, P.~Moree and R.~Osburn, Principal forms $X^2+nY^2$ representing many integers,
\textit{Abh.~Math.~Sem.~Univ.~Hambg.}\, \textbf{81} (2011), 129--139.

\bibitem{CMM} S.\,Charlton, A.\,Medvedovsky and P.\,Moree, Euler--Kronecker constants of modular forms:\,beyond
Dirichlet $L$-series, \url{https://arxiv.org/abs/2412.01803}, 
submitted for publication.


\bibitem{C} J.~Cilleruelo, The least common multiple of a quadratic sequence, \textit{Compos. Math.}
\textbf{147} (2011),
1129--1150.

\bibitem{CRSZ} J.~Cilleruelo, J.~Rué, P.~\v{S}arka and A.~Zumalac\'arregui, 
The least common multiple of random sets of positive integers, 
\textit{J. Number Theory} \textbf{144} (2014), 92--104.

\bibitem{CLM} A.~Ciolan, A.~Languasco and P.~Moree, Landau and Ramanujan approximations for divisor sums and coefficients of cusp forms,  \textit{J.\,Math.\,Anal.\,Appl.} \textbf{519} (2023), no. 2, Paper No. 126854.



\bibitem{CK1} R.D.~Connors and J.P.~Keating,
Two-point spectral correlations for the square billiard,
\textit{J. Phys. A} \textbf{30} (1997), 1817--1830.

\bibitem{CK2} R.D.~Connors and J.P.~Keating, Degeneracy moments for the square billiard,
\textit{J. Phys. A} \textbf{32} (1999), 555--562.


\bibitem{conway} J.H.~Conway and N.J.A.~Sloane, Lattices with
few distances, \textit{J.~Number Theory} \textbf{39} (1991), 75--90.

\bibitem{Daas} \pieter{M.~Daas, Congruence conditions for the mod $\lambda$ vanishing of the Fourier coefficients of classical eigenforms, 
\url{https://arxiv.org/abs/2506.08865}, submitted for publication.}

\bibitem{DG}
S.~Das and S.~Ganguly, Gaps between non-zero Fourier coefficients of cusp forms, \textit{Proc. Amer. Math. Soc.} 
\textbf{142} (2014), 3747--3755.

\bibitem{DDNS} C.~David, L.~Devin, J.~Nam and J.~Schlitt,
Lemke Oliver and Soundararajan bias for consecutive sums of two squares,
\textit{Math. Ann.} \textbf{384} (2022), 1181--1242.

\bibitem{BT3}
R.~de la Bret\`eche and G.~Tenenbaum, 
Remarks on the Selberg--Delange method,
\textit{Acta Arith.} \textbf{200} (2021), 349--369.


\bibitem{BT}
R.~de la Bret\`eche and G.~Tenenbaum, 
Friable averages of oscillating arithmetic functions, \textit{Number theory in memory of Eduard Wirsing}, 43--71, Springer, Cham, 2023.

\bibitem{BT2}
R.~de la Bret\`eche and G.~Tenenbaum,
Friable averages of complex arithmetic functions,
\textit{Funct. Approx. Comment. Math.} \textbf{71} (2024), 87--119.



\bibitem{DelShi} P.~Deligne, Travaux de Shimura, 
S\'eminaire Bourbaki, 23\`eme ann\'ee (1970/1971), Exp. No. 389, pp. 123--165,
\textit{Lecture Notes in Math.} \textbf{244}, Springer, Berlin-New York, 1971.

\bibitem{Deligne} P.~Deligne, La conjecture de Weil. I,
\textit{Inst. Hautes \'Etudes Sci. Publ. Math.} \textbf{43} (1974), 273--307.

\bibitem{DHZ} M.~Derickx, M.~van Hoeij and J.~Zeng,
Computing Galois representations and equations for modular
curves $X_H(\ell),$
\url{https://arxiv.org/pdf/1312.6819.pdf}, unpublished
preprint, 2013.



\bibitem{DS}
F.~Diamond and J.~Shurman, \textit{A First Course in Modular Forms},
Graduate Texts in Mathematics \textbf{228}, Springer-Verlag, New York, 2005. 

\bibitem{Dyson} F.J.~Dyson, A walk through Ramanujan's garden, \textit{Ramanujan Revisited} (Urbana-Champaign, Ill., 1987), Academic Press, Boston, MA, 1988, 7--28.


\bibitem{taucomp}  B.~Edixhoven and J.-M.~Couveignes (Eds.),
\textit{Computational aspects of modular forms and Galois representations.
How one can compute in polynomial time the value of Ramanujan's tau at a prime},  Annals of Mathematics Studies
\textbf{176}, Princeton University Press, Princeton, NJ, 2011.


\bibitem{EGS}
P.~Erd\H{o}s, A.~Granville 
and C.~Spiro, On the normal behavior of the iterates of some arithmetic functions. \textit{Analytic number theory} (Allerton Park, IL, 1989), 165--204, \textit{Progr. Math.} \textbf{85}, Birkh\"auser Boston, Boston, MA, 1990.

\bibitem{ES}
P.~Erd\H{o}s and S.~Szekeres, \"Uber die Anzahl der Abelchen Gruppen gegebener Ordnung und \"uber ein verwandtes zahlentheorestisches Problem, 
\textit{Ada Scient. Math. Szeged} \textbf{7} (1935), 95--102.

\bibitem{Estermann}
T.~Estermann, An asymptotic formula in the theory of numbers,
\textit{Proc. London Math. Soc.} (2) \textbf{34} (1932), 280--292.

\bibitem{ERS} S.~Ettahri, O.~Ramar\'e and L.~Surel, Fast multi-precision computation of some
Euler products, \textit{Math. Comp.} \textbf{90} (2021),
2247--2265.

\bibitem{constants} S.R.~Finch, \textit{Mathematical constants}, Encyclopedia of Mathematics and its Applications \textbf{94}, Cambridge University Press, Cambridge, 2003.

\bibitem{FMS} S.R.~Finch, G.~Martin and P.~Sebah, Roots of unity and nullity
modulo $n$, \textit{Proc. Amer. Math. Soc.} \textbf{138}  (2010), 2729--2743.

\bibitem{Fomenko} O.M.~Fomenko, Distribution of 
values of 
Fourier coefficients of modular forms of weight 1, 
\textit{J. Math. Sci. (New York)} \textbf{89} 
(1998), 1050--1071.

\bibitem{FLM} K.~Ford, F.~Luca and P.~Moree, Values of the Euler $\phi$-function not divisible by
a given odd prime, and the distribution of Euler--Kronecker constants for
cyclotomic fields, \textit{Math. Comp.} \textbf{83} (2014), 1447--1476.

\bibitem{FK} \'E.~Fouvry and J.~Kl\"uners, 
On the negative Pell equation, 
\textit{Ann. of Math.} (2) \textbf{172} (2010), 2035--2104.

\bibitem{FLW} \'E.~Fouvry, C.~Levesque
 and M.~Waldschmidt, Representation of integers by cyclotomic binary 
 forms, \textit{Acta Arith.} \textbf{184} (2018), 67--86.

\bibitem{FKR} T.~Freiberg, P.~Kurlberg and L.~Rosenzweig,
Poisson distribution for gaps between sums of two squares and level spacings for toral point scatterers,
\textit{Commun. Number Theory Phys.} \textbf{11} (2017),
837--877.


\bibitem{Gaussconstant} Gauss's constant, \url{https://mathworld.wolfram.com/GausssConstant.html}.



\bibitem{Gelbart} S.~Gelbart,
An elementary introduction to the Langlands program,
\textit{Bull. Amer. Math. Soc. (N.S.)} \textbf{10} (1984), 177--219.



\bibitem{poly} O.~Gorodetsky, A polynomial analogue of Landau's theorem and related problems,
\textit{Mathematika} \textbf{63}
(2017), 622--665.

\bibitem{Goro2} O.~Gorodetsky, Sums of two squares are strongly
biased towards quadratic residues,
\textit{Algebra Number Theory} \textbf{17}
(2023), 775--804.

\bibitem{GoroRodgers} O.~Gorodetsky and B.~Rodgers, 
The variance of the number of sums of two squares in 
$\mathbb F_q[T]$ in short intervals, 
\textit{Amer. J. Math.} \textbf{143} (2021), 1703--1745.

\bibitem{GK} A.~Granville and D.~Koukoulopoulos, Beyond the LSD method for the
partial sums of multiplicative functions,
\textit{Ramanujan J.} \textbf{49} (2019), 287--319.

\bibitem{GS} \pieter{A.~Granville and 
K.~Soundararajan, An uncertainty principle for arithmetic sequences,
\textit{Ann. of Math.} (2) \textbf{165} 
(2007), 593--635.}

\bibitem{Gray}
J.J.~Gray, A commentary on Gauss's mathematical diary, 1796--1814, with an English translation,
\textit{Exposition. Math.} \textbf{2} (1984), 97--130.

\bibitem{Grosswald} E.~Grosswald, \textit{Representations of Integers as Sums of Squares}, Springer-Verlag, New York, 1985.

\bibitem{HTW} \pieter{G.~Hanrot, G.~Tenenbaum and J.~Wu,
Moyennes de certaines fonctions multiplicatives sur les entiers friables. II,
\textit{Proc. Lond. Math. Soc.} (3) \textbf{96} (2008), 107--135.}

\bibitem{hardysquares}
G.H.~Hardy, On the expression of a number as the sum of two
squares, \textit{Quart.~J.~Math.~(Oxford)} \textbf{46} (1915), 263--283.


\bibitem {Hardy} G.H.~Hardy, \textit{Ramanujan}, 3rd ed., Chelsea, New York,
1978.

\bibitem{HL} G.~Hardy and J.E.~Littlewood, Some problems of `Partitio numerorum'; III: On the expression of a number as a sum of primes, \textit{Acta Math.} \textbf{44} (1923), 1--70.

\bibitem{hw}
G.H.~Hardy and E.M.~Wright, \textit{An Introduction to the Theory of Numbers}, 6th.~ed., Oxford Univ.~Press, 2008.

\bibitem {Hare} D.E.G.~Hare, The
Landau--Ramanujan constant $K$ (with 125.079  decimals), \url{https://oeis.org/A064533/a064533_1.txt}.

\bibitem{Havil} G.~Havil,
\textit{Gamma. Exploring Euler's Constant},
Princeton University Press, Princeton, NJ, 2003.

\bibitem{H-B} D.R.~Heath-Brown, Ternary quadratic forms and sums of three square-full 
numbers, \textit{S\'eminaire de Th\'eorie des Nombres}, Paris 1986--87, ed.\,C.~Goldstein (Birkh\"auser, Boston, 1988), 137--163.

\bibitem{H-B2} D.R.~Heath-Brown,  
Equidistribution for solutions of $p+m^2+n^2=N$, and for 
Ch\^atelet surfaces, \textit{Acta Arith.} \textbf{214} (2024), 
23--37.

\bibitem{Hecke} E.~Hecke, \"Uber Modulfunktionen und die Dirichletschen Reihen mit Eulerscher Produktentwicklung. I,
\textit{Math. Ann.} \textbf{114} (1937),
1--28.

\bibitem {Heupel} W.~Heupel, Die Verteilung der ganzen Zahlen, die
durch quadratische Formen dargestellt werden, \textit{Arch.~Math.}~\textbf{19}
(1968), 162--166.

\bibitem{HQT} S.~Hong, G.~Qian and Q.~Tan, The least common multiple of a sequence of products of linear polynomials, 
\textit{Acta Math. Hungar.} \textbf{135} 
(2012), 160--167.




\bibitem{Hooley57} C.~Hooley,
On the representation of a number as the sum of two squares and a prime, \textit{Acta Math.} \textbf{97} (1957), 189--210.

\bibitem{Hooley1} C.~Hooley,
On the intervals between numbers that are sums of two squares,
\textit{Acta Math.} \textbf{127} (1971), 279--297.

\bibitem{Hooley2} C.~Hooley,
On the intervals between numbers that are sums of two squares. II, 
\textit{J. Number Theory} 
\textbf{5} (1973), 215--217.

\bibitem{Hooley3} C.~Hooley,
On the intervals between numbers that are sums of two squares. III.
\textit{J. Reine Angew. Math.} \textbf{267} (1974), 207--218.

\bibitem{Hooleybook} C.~Hooley, 
\textit{Applications of Sieve Methods to the Theory of Numbers}, Cambridge Tracts in Mathematics \textbf{70}, Cambridge Univ. Press, Cambridge-New York-Melbourne, 1976.

\bibitem{Hooley4} C.~Hooley, On the intervals between numbers that are sums of two squares. IV,
\textit{J. Reine Angew. Math.}
\textbf{452} (1994), 79--109.


\bibitem{huxley}
M.N.~Huxley, Exponential sums and lattice points. III, \emph{Proc.~London
Math.~Soc}.~(3) \textbf{87} (2003), 591--609.


\bibitem{HIS} D.~Hu, H.~Iyer and A.~Shashkov, Modular forms and an explicit Chebotarev variant of the Brun-Titchmarsh theorem,
\textit{Res. Number Theory} \textbf{9}, Paper No. 46, 37 p. (2023).

\bibitem{Ihara} Y.~Ihara,  
On the Euler--Kronecker constants of global fields and primes with small norms, \textit{Algebraic geometry and number theory}, 407--451, \textit{Progr. Math.} \textbf{253},
Birkh\"auser Boston, Inc., Boston, MA, 2006.

\bibitem{Indlekofer}
K.-H.~Indlekofer, Scharfe untere Absch\"atzung f\"ur die 
Anzahlfunktion der B-Zwillinge, \textit{Acta Arith.} 
\textbf{26} (1974/75), 207--212.

\bibitem {Iwaniec} H.~Iwaniec, The half dimensional sieve, \textit{Acta Arith.} \textbf{ 29} (1976), 69--95.

\bibitem{jacobi}
C.G.J.~Jacobi, \textit{Fundamenta Nova Theoriae Functionum Ellipticarum}, Sumptibus Fratrum Borntr\"{a}ger, Regiomonti, 1829. 

\bibitem {James} R.D.~James, The distribution of integers represented
by quadratic forms, \textit{Amer. J. Math.} \textbf{60} (1938), 737--744.

\bibitem{Joyner} D.~Joyner, 
\textit{Distribution theorems of 
$L$-functions},
Pitman Research Notes in Mathematics Series 
\textbf{142}, Longman Scientific \& Technical, Harlow; John Wiley \& Sons, Inc., New York, 1986. 

\bibitem{GSS} A. Granville, A. Sedunova and C. Sabuncu, 
The multiplication table constant and sums of two squares, 
\textit{Acta Arith.} \textbf{214} (2024), 499--522.

\bibitem{Kanigel}
R.~Kanigel, \textit{The man who knew infinity.
A life of the genius Ramanujan}, Charles Scribner's Sons, New York, 1991.

\bibitem{Katz} N.M.~Katz, An overview of Deligne's proof of the Riemann hypothesis for varieties over finite fields, \textit{Mathematical developments arising from Hilbert problems} (Proc. Sympos. Pure Math., Vol. XXVIII, Northern Illinois Univ., De Kalb, Ill., 1974), 275--305, Amer. Math. Soc., Providence, R.I., 1976.



\bibitem{KK} N.~Kimmel and V.~Kuperberg, Positive density for consecutive runs of sums of two 
squares, \url{https://arxiv.org/pdf/2406.04174}.

\bibitem{Koblitz} N.~Koblitz, \textit{Introduction to Elliptic Curves and Modular Forms}, Second edition, Graduate Texts in Mathematics \textbf{97}, Springer-Verlag, New York, 1993. 

\bibitem{Koukou} D.~Koukoulopoulos, \textit{The distribution of prime numbers},
Graduate Studies in Mathematics \textbf{203},
American Mathematical Society, Providence, RI, 2019.



\bibitem{Kowalski} E.~Kowalski, Some aspects and applications of the Riemann hypothesis over finite fields, \textit{Milan J. Math.} \textbf{78}
(2010), 179--220.

\bibitem{KP} P.~Koymans and C.~Pagano, On Stevenhagen's conjecture,
\url{https://arxiv.org/abs/2201.13424}, submitted for publication.

\bibitem{Lagarias} J.C.~Lagarias, Euler's constant: Euler's work and modern developments,
\textit{Bull. Amer. Math. Soc. (N.S.)} \textbf{50} (2013), 527--628.

\bibitem {L} E.~Landau, \"Uber die Einteilung
der positiven ganzen Zahlen in vier Klassen nach der mindest Anzahl
der zu ihrer additiven Zusammensetzung erforderlichen Quadrate,
\textit{Arch. der Math. und Phys.} (3) \textbf{13} (1908),
305--312; (see also \textit{Collected Works} \textbf{4}, 59--66).

\bibitem{La09a} E.~Landau, L\"osung des Lehmer'schen Problems,
\textit{Amer.~J.~Math.} \textbf{31} (1909), 86--102;
(see also \textit{Collected Works} \textbf{4}, 131--147).

\bibitem {Lehrbuch} E.~Landau,  \textit{Handbuch der Lehre von der
Verteilung der Primzahlen.} II, Teubner, Leipzig, 1909; 2nd ed., Chelsea,
New York, 1953.

\bibitem {Langlands} R.P.~Langlands, \textit{Euler Products},
A James K. Whittemore Lecture in Mathematics given at Yale University, 1967, 
\textit{Yale Mathematical Monographs} 
\textbf{1}, Yale University Press, New Haven, 
Conn.-London, 1971. 

\bibitem {Sconstant} 
A.~Languasco, Second-order term in the asymptotic expansion of $B(x)$, the count of numbers up to $x$ which are the sum of two squares (with 130,000 decimals), \url{https://oeis.org/A227158},
see also \url{https://www.dei.unipd.it/~languasco/CLM.html}.


\bibitem{Languasco2021a}
A.~Languasco, Efficient computation of the {E}uler-{K}ronecker constants of prime
  cyclotomic fields,
  \textit{Res. Number Theory}
  \textbf{7} (2021), no.\,1, Paper No.\,2, 22 pp.

\bibitem{Languasco2023a}
A.~Languasco,
A unified strategy to compute some special functions of number-theoretic interest.
\textit{J. Number Theory} \textbf{247} (2023), 118--161.

\bibitem{LaMoree} \pieter{A.~Languasco and P.~Moree,
Euler constants from primes in arithmetic progression, 
\textit{Math.\,Comp.} (2025), \url{https://doi.org/10.1090/mcom/4057}.}


\bibitem{LanguascoR2021}
A.~Languasco and L.~Righi,
A fast algorithm to compute the {R}amanujan-{D}eninger gamma function
  and some number-theoretic applications.
\textit{Math. Comp.} \textbf{90} (2021), 2899--2921.



\bibitem{Lehmertable} D.H.~Lehmer,
Ramanujan's function $\tau(n),$ \textit{Duke Math. J.}
\textbf{10} (1943), 483--492.


\bibitem{Lehmer} D.H.~Lehmer, The vanishing of
Ramanujan's function $\tau(n),$ \textit{Duke Math. J.}
\textbf{14} (1947), 429--433.

\bibitem{LOSound} R.J.~Lemke Oliver and
K.~Soundararajan, Unexpected biases in the distribution of consecutive primes,
\textit{Proc. Natl. Acad. Sci. USA} \textbf{113} (2016), no. 31, E4446--E4454.

\bibitem{LeVeque} W.J.~LeVeque, 
The distribution of values of multiplicative functions, 
\textit{Michigan Math. J.} \textbf{2} (1953/54), 179--192 (1955).

\bibitem{Levy} P.~L\'evy, Observations sur le m\'emoire de M.F. Tricomi, \textit{Atti Accad. Sci. Torino Cl. Sci. Fis. Mat. Nat.} \textbf{75} (1939), 177--183.

\bibitem{winnie} W.-C.W.~Li, The Ramanujan conjecture and its applications, \textit{Philos. Trans. Roy. Soc. A} \textbf{378} (2020), no. 2163, 20180441, 14 pp.

\bibitem{xx}
X.~Li and X.~Yang, An Improvement on Gauss’s Circle Problem and Dirichlet’s Divisor Problem, to appear.

\bibitem{Linnik} Yu.V.~Linnik,
All large numbers are sums of a prime and two squares (A problem of Hardy and Littlewood). I, 
\textit{Mat. Sb. (N.S.)} 52(94) (1960), 661--700; II.~ibid.~53(95) (1961), 3--38.



\bibitem{LR} N.~Lygeros and O.~Rozier,  A new solution to the equation $\tau(p)\equiv 0(\text{mod~}p)$,
\textit{J. Integer Seq.} \textbf{13} (2010), no. 7, Article 10.7.4, 11 pp.

\bibitem{Maier} \pieter{H.~Maier, 
Chains of large gaps between consecutive primes, 
\textit{Adv. in Math.} \textbf{39} 
(1981), 257--269}.

\bibitem{MR} J.~Maynard and Z.~Rudnick, A lower bound on the least common multiple of polynomial sequences, 
\textit{Riv. Math. Univ. Parma (N.S.)} 
\textbf{12} (2021), 143--150.

\bibitem{error} B.~Mazur,
Finding meaning in error terms,
\textit{Bull. Amer. Math. Soc.} (N.S.) \textbf{45} (2008), 185--228.

\bibitem{Milne} S.C.~Milne, Infinite families of exact sums of squares formulas, Jacobi elliptic functions, continued fractions, and Schur functions,  \textit{Ramanujan J.} \textbf{6} (2002), 7--149; published under the same title in Dev.~Math., No.~5, Kluwer Academic Publishers, Boston, MA, 2002, 149 pp.

\bibitem{MRama} P.~Moree, On some claims in Ramanujan's `unpublished'
manuscript on the partition and tau functions,
\textit{Ramanujan J.} \textbf{8} (2004), 317--330.

\bibitem{Mpreprint} P. Moree, Values of the Euler phi function not divisible by a prescribed odd prime,
unpublished preprint (precursor of \cite{FLM}) available at \url{https://arxiv.org/abs/math/0611509}.

\bibitem{LvR} P.~Moree, Counting numbers in multiplicative sets:\ Landau versus Ramanujan, 
\textit{Math. Newsl.} \#3 (2011), 73--81.


\bibitem{MC} P.~Moree and J.~Cazaran, On a claim of Ramanujan in
his first letter to Hardy, \textit{Expos. Math.} \textbf{17} (1999),
289--312.

\bibitem{MN} \pieter{P.~Moree and A.~Noubissie, Higher reciprocity laws 
and ternary linear recurrence sequences,
\textit{Int. J. Number Theory} \textbf{21} 
 (2025), 993--1015.}

\bibitem{MO} P.~Moree and R.~Osburn, Two-dimensional lattices
with few distances, \textit{Enseignement Math.} \textbf{52} (2006), 361--380.

\bibitem{MteR} P.~Moree and H.J.J.~te Riele, The hexagonal versus
the square lattice, \textit{Math.~Comp.} \textbf{73} (2004), 451--473.

\bibitem{MurtyRR} M.R.~Murty, The Ramanujan $\tau$ function, in
\textit{Ramanujan revisited} (Urbana-Champaign, Ill., 1987), 269--288, Academic Press, Boston, MA, 1988.

\bibitem{Murtysieving} M.R.~Murty, 
Sieving using Dirichlet series, 
\textit{Currents trends in number theory (Allahabad, 2000)}, 
111--124, Hindustan Book Agency, New Delhi, 2002.

\bibitem{MM} M.R.~Murty and V.K.~Murty,
\textit{The Mathematical Legacy of Srinivasa Ramanujan}, Springer, New Delhi, 2013.

\bibitem{MMS} M.R.~Murty, V.K.~Murty and
T.N. Shorey, Odd values of the Ramanujan $\tau$-function, \textit{Bulletin Soc.
Math. France} \textbf{115} (1987),
391--395.

\bibitem {O4} R.W.K.~Odoni, On norms of integers in a full module of
an algebraic number field and the distribution of values of binary integral
quadratic forms, \textit{Mathematika} \textbf{22} (1975), 108--111.

\bibitem{Odonishort} R.W.K.~Odoni,
Solution of some problems of Serre on modular forms: The method of Frobenian functions, 
\textit{Recent progress in analytic number theory}, Symp. Durham 1979, Vol. 2 (1981), 159-169 (1981).


\bibitem {O1} R.W.K.~Odoni, Notes on the method of Frobenian functions with
applications to Fourier coefficients of modular forms, \textit{Banach Center Publications} \textbf{17},
PWN, Warsaw, 1985.




\bibitem{O} R.W.K.~Odoni, A problem of Rankin on sums of powers
of cusp-form coefficients, \textit{J.~London Math.~Soc.}\, \textbf{44}
(1991), 203--217.

\bibitem{OEIS} OEIS, Padovan sequence, 
\url{https://oeis.org/A000931}.

\bibitem{Onoweb} K.~Ono, \textit{The Web of Modularity: Arithmetic of the Coefficients of 
Modular Forms and $q$-series},
CBMS Regional Conference Series in Mathematics \textbf{102}, Published for the Conference Board of the Mathematical Sciences, Washington, DC; by the American Mathematical Society, Providence, RI, 2004.



\bibitem{Pall2} G.~Pall, The distribution of integers represented by
binary quadratic forms, \textit{Bull.~Am.~Math.~Soc.}\, \textbf{49} (1943), 447--449.

\bibitem{PS}
T.R.~Parkin and  D.~Shanks, On the distribution of parity in the partition 
function, \textit{Math. Comp.} \textbf{21} (1967), 466--480.

\bibitem {P} A.G.~Postnikov, \textit{Introduction to Analytic Number
Theory}, AMS translations of mathematical monographs \textbf{68}, AMS,
Providence, Rhode Island, 1988.

\bibitem{Prachar} K.~Prachar, \"Uber Zahlen der Form $a^2+b^2$ in einer arithmetischen Progression, \textit{Math.~Nachr.}\, \textbf{10} (1953), 51--54.

\bibitem{Prachar2} K.~Prachar,
Bemerkung zu der Arbeit von Herrn H. Beki\'c,
\textit{J. Reine Angew. Math.} \textbf{229} (1968), 28.

\bibitem{Radu}
S.~Radu, A proof of Subbarao’s conjecture, 
\textit{J. Reine Angew. Math.} \textbf{672} (2012), 161--175.

\bibitem{Ramachandra} 
K.~Ramachandra, 
Some problems of analytic number theory, 
\textit{Acta Arith.} \textbf{31} (1976), 313--324.

\bibitem{Rama1916} S.~Ramanujan, On certain arithmetical
functions, \textit{Trans. Cambridge Phil. Soc.} \textbf{22}
(1916), 159--184.

\bibitem{nb}
S.~Ramanujan, \textit{Notebooks} (2 volumes), Tata Institute of Fundamental
Research, Bombay, 1957; second ed., 2012.

\bibitem{cp}
S.~Ramanujan, \textit{Collected Papers}, Cambridge Univ.~Press, Cambridge, 1927; reprinted by Chelsea, New York, 1962; reprinted by the American Mathematical Society, Providence, 2000.

\bibitem{lnb}
S.~Ramanujan, \textit{The Lost Notebook and Other Unpublished
Papers}, Narosa, New Delhi, 1988.

\bibitem{Accurate} O.~Ramar\'e, Accurate computations of Euler products over primes in arithmetic progressions,
\textit{Funct. Approx. Comment.~Math.} \textbf{65} (2021), 33--45.

\bibitem{walk} O.~Ramar\'e,
\textit{Excursions in Multiplicative Number Theory} (with contributions
by P.~Moree and A.~Sedunova),
Birkh\"auser Advanced Texts: Basler Lehrb\"ucher, Birkh\"auser/Springer, Basel, 2022.



\bibitem{Rankin61} R.A.~Rankin, The divisibility of divisor functions,
\textit{Proc.~Glasgow Math.~Assoc.}\, \textbf{5} (1961), 35--40.

\bibitem{watsonobi} R.A.~Rankin,
George Neville Watson (obituary),
\textit{J.~London Math.~Soc.}\,
\textbf{41} (1966), 551--565.

\bibitem{Rankin76} R.A.~Rankin,
Ramanujan's unpublished work on congruences,
\textit{Modular functions of one variable.~V} (Proc.~Second Internat.~Conf., Univ.~Bonn, Bonn, 1976), pp. 3–15,
\textit{Lecture Notes in Math.} \textbf{601}, Springer, Berlin, 1977.



\bibitem{rankinzelf} R.A.~Rankin, Ramanujan's manuscripts and notebooks,
{\it Bull.~London Math.~Soc.} \textbf{14} (1982), 81--97.

\bibitem{rankintau} R.A.~Rankin,
Ramanujan's tau-function and its generalizations,
\textit{Ramanujan revisited} (Urbana-Champaign, Ill., 1987), Academic Press, Boston, MA, 1988, 245--268.

\bibitem{Rieger} G.J.~Rieger, Aufeinanderfolgende Zahlen als Summen von zwei Quadraten, 
Nederl. Akad. Wetensch. Proc. Ser. A 68,
\textit{Indag. Math.} \textbf{27} (1965), 208--220.


\bibitem{Robin} G.~Robin,
Grandes valeurs de la fonction somme des diviseurs et hypoth\`ese de Riemann,
\textit{J. Math. Pures Appl.} (9)
\textbf{63} (1984), 187--213.


\bibitem{RS} J.B.~Rosser and L.~Schoenfeld, Approximate formulas for some functions of
prime numbers, \textit{Illinois Journal Math.} \textbf{6} (1962), 64--94.

\bibitem{RZ} Z.~Rudnick and S.~Zehavi, On Cilleruelo's conjecture for the least common multiple of polynomial sequences, 
\textit{Rev. Mat. Iberoam.} \textbf{37} (2021), 1441--1458.

\bibitem{Sabuncu}
C.~Sabuncu, On the moments of the number of representations as sums of two prime squares, \textit{Int. Math. Res. Not. IMRN 2024} \textbf{11}, 9411--9439.

\bibitem{Schlitt}
J.~Schlitt, Biases towards the zero residue class for quadratic forms in arithmetic progressions, \url{https://arxiv.org/abs/2308.13959}.

\bibitem {Schmutz} P.~Schmutz Schaller, Geometry of Riemann surfaces
based on closed geodesics, \textit{Bull.~Am.~Math.~Soc.}\, \textbf{35} (1998),
193--214.



\bibitem{Scourfield64} E.J.~Scourfield, On the divisibility of $\sigma_{\nu}(n)$, \textit{Acta Arith.}\,  \textbf{10} (1964), 245--285.


\bibitem{Scourfield72} E.J.~Scourfield, 
Non-divisibility of some multiplicative functions, 
\textit{Acta Arith.} \textbf{22} (1972/73), 287--314.

\bibitem{Scourfield77} E.J.~Scourfield, 
On the divisibility of $r_k(n)$, 
\textit{Glasgow Math. J.} \textbf{18} (1977), 109--111.

\bibitem{Sedunova22} A.~Sedunova, Intersections of binary quadratic forms in primes and the paucity phenomenon, 
\textit{J. Number Theory} \textbf{235} (2022), 305--327.

\bibitem{Selberg} \pieter{A.~Selberg, \textit{Collected Papers}, Vol. II, Springer-Verlag, Berlin, 1991.}

\bibitem{Serre73} J-P.~Serre,
Congruences et formes modulaires [d'apr\`es H.P.F.~Swinnerton-Dyer],
S\'eminaire Bourbaki, 24\`eme ann\'ee (1971/1972), Exp. No. 416, pp. 319--338,
\textit{Lecture Notes in Math.} \textbf{317}, Springer, Berlin-New York, 1973.

\bibitem{course} J-P.~Serre, \textit{A Course in Arithmetic},
Graduate Texts in Mathematics \textbf{7}, Springer-Verlag, New York-Heidelberg, 1973.


\bibitem {Ser} J-P.~Serre, Divisibilit\'e de certaines
fonctions arithm\'etiques, \textit{Enseignement Math.}\, \textbf{22} (1976),
227--260.



\bibitem{SerreChebotarev} J-P.~Serre, Quelques applications du th\'eor\`eme de densit\'e de Chebotarev, \textit{Inst. Hautes
\'Etudes Sci. Publ. Math.} \textbf{54} (1981), 323--401.

\bibitem{serrelacunary} 
J-P.~Serre, Sur la lacunarit\'e des puissances de 
$\eta$, \textit{Glasgow Math. J.} \textbf{27} (1985), 203--221. 

\bibitem{SerreJordan} J-P.~Serre,
On a theorem of Jordan,
\textit{Bull. Amer. Math. Soc. (N.S.)} \textbf{40} (2003), 429--440.

\bibitem{Sah} A.~Sah, An improved bound on the least common multiple of polynomial sequences, \textit{J. Th\'eor. Nombres Bordeaux} 
\textbf{32} (2020), 891--899.

\bibitem{Shanks} D.~Shanks, The second-order term in the asymptotic
expansion of $B(x)$, \textit{Math.\,Comp.}\, \textbf{18} (1964), 75--86.

\bibitem{SS} D.~Shanks and L.P.~Schmid, Variations on a theorem
of Landau, I,
{\it Math.~Comp.}\, \textbf{20} (1966), 551--569.

\bibitem {Shiu}  P.~Shiu, Counting sums of two squares: the
Meissel-Lehmer method, \textit{Math.~Comp.}\, \textbf{47} (1986), 351-360: Corrigendum: \textit{Math.~Comp.}\, \textbf{88} (2019), 2935--2938.

\bibitem{siegel}
C.L.~Siegel, \textit{Topics in Complex Function Theory}, Vol.~1, Wiley, New York, 1969.

\bibitem{Song}
J.M.~Song, Sums of nonnegative multiplicative functions over integers without large prime factors. I,
\textit{Acta Arith.} 
\textbf{97} (2001), 329--351; 
II, ibid, \textbf{102} (2002), 105--129.


\bibitem {Stanley} G.K.~Stanley, Two assertions made by Ramanujan,
\textit{J.~London Math.~Soc.}\, \textbf{3} (1928), 232--237. Corrigenda, {\it ibid}
\textbf{4} (1929), 32.

\bibitem {peterpell} P.\,Stevenhagen, The number of real quadratic fields
having units of negative norm, \textit{Experiment. Math.} 
\textbf{2} (1993), 121--136.


\bibitem{SL} P.~Stevenhagen and H.W. Lenstra, Jr., 
Chebotar\"ev and his density theorem, 
\textit{Math. Intelligencer} \textbf{18} (1996), 26--37.

\bibitem{SX}
C.L.~Stewart and Y.S.~Xiao, On the representation of integers by binary forms,
\textit{Math.~Ann.}\, \textbf{375} (2019), 133--163.


\bibitem{subbarao}
M.~Subbarao, Some remarks on the partition function, 
\textit{Amer.~Math.~ Monthly} \textbf{73} (1966), 851–-854.

\bibitem{S-D-ladic}  H.P.F.~Swinnerton-Dyer,
On $\ell$-adic representations and congruences for coefficients of modular forms, 
\textit{Modular functions of one variable} III (Proc. Internat. Summer School, Univ. Antwerp, 1972), \textit{Lecture Notes in Math.} 
\textbf{350}, Springer, Berlin, 1973, 1--55.


\bibitem{S-D} H.P.F.~Swinnerton-Dyer, Congruence properties of $ \tau(n),$ \textit{Ramanujan Revisited} (Urbana-Champaign, Ill., 1987), 289--311, Academic Press, Boston, MA, 1988.

\bibitem{Tenenbaum} G.~Tenenbaum, \textit{Introduction to Analytic and Probabilistic Number Theory}, third edition, Graduate Studies in Mathematics
\textbf{163}, American Mathematical Society, Providence, RI, 2015.

\bibitem{Tenenbaumeffective} G.~Tenenbaum,  
Moyennes effectives de fonctions multiplicatives complexes, \textit{Ramanujan J.} \textbf{44} (2017), 641--701.

\bibitem{TW} \pieter{G.~Tenenbaum and J.~Wu,
Moyennes de certaines fonctions multiplicatives sur les entiers friables,
\textit{J. Reine Angew. Math.} \textbf{564} (2003), 119--166.}

\bibitem{Troupe} L.~Troupe, Divisor sums representable as the sum of two squares, 
\textit{Proc. Amer. Math. Soc.} \textbf{148} (2020), 4189--4202.



\bibitem{Blij} F.~van der Blij, Binary quadratic forms of discriminant $-23$,
\textit{Nederl. Akad. Wetensch. Proc. Ser.} A \textbf{55},
\textit{Indag. Math.} \textbf{14} (1952), 498--503.

\bibitem{gnwatson}
G.N.~Watson, Scraps from some mathematical note books,
\textit{Math.~Gazette} \textbf{18} (1934), 5--18.

\bibitem{watson1} G.N.~Watson, \"Uber Ramanujansche Kongruenzeigenschaften der
Zerf\"allungsanzahlen. I, \textit{Math.~Z.}\, \textbf{39} (1935), 712--731.



\bibitem{watson3} G.N.~Watson, A table of Ramanujan's function $\tau(n)$,
\textit{Proc. London Math. Soc.} \textbf{51} (1949), 1--13.


\bibitem {Williams2} K.S.~Williams, Note on integers representable
by binary quadratic forms, \textit{Canad.~Math.~Bull.}\, \textbf{18} (1975), 123--125.

\bibitem{Wilton} J.R.~Wilton,
Congruence properties of Ramanujan's function $\tau(n)$,
\textit{Proc. London Math. Soc.} (2) \textbf{31} (1930), 1--10.


\bibitem {Wirsing} E.~Wirsing, Das asymptotische Verhalten von Summen
\"uber multiplikative Funktionen, \textit{Math.~Ann.}\, \textbf{143} (1961),
75--102.



\bibitem {Wirsing2} E.~Wirsing,
Das asymptotische Verhalten von Summen \"uber 
multiplikative Funktionen. II, 
\textit{Acta Math. Acad. Sci. Hung.} \textbf{18} 
(1967), 411--467.

\end{thebibliography}
\end{document}